\documentclass[onefignum,onetabnum]{siamart190516}
\usepackage{amsmath,amssymb,amsxtra}
\usepackage{mathrsfs}
\usepackage{float,verbatim}
\usepackage{threeparttable,booktabs}
\usepackage{graphicx}
\usepackage{epstopdf}
\usepackage{mathtools}
\usepackage{subfig}
\usepackage{algorithm,algorithmic}
\usepackage{color}
\graphicspath{{./pics/}}

\usepackage{array}
\newtheorem{remark}{Remark}[section]

\usepackage{graphicx}
\usepackage{dsfont}
\usepackage{epstopdf}
\graphicspath{{./pics/}}
\usepackage{booktabs}
\usepackage{threeparttable}
\usepackage{verbatim}
\newcommand{\bel}{\begin{equation} \label}
\newcommand{\ee}{\end{equation}}
\def\beq{\begin{equation}}
\def\eeq{\end{equation}}
\newcommand{\bea}{\begin{eqnarray}}
\newcommand{\eea}{\end{eqnarray}}
\newcommand{\beas}{\begin{eqnarray*}}
\newcommand{\eeas}{\end{eqnarray*}}

\newcommand{\p}{\partial}

\newcommand{\eqnref}[1]{(\ref {#1})}

\allowdisplaybreaks

\newtheorem{thm}{Theorem}[section]

\newtheorem{assumption}{Assumption}[section]

\newtheorem{cor}{Corollary}[section]
\numberwithin{equation}{section}

\allowdisplaybreaks

\def\phi {\varphi}

\allowdisplaybreaks

\allowdisplaybreaks

\title{Identification of a Spatially-Dependent Variable Order in One-Dimensional Subdiffusion\thanks{The work of B. Jin is supported by Hong Kong RGC General Research Fund (Project
14306423), and a start-up fund from The Chinese University of Hong Kong.}}

\author{Jiho Hong\thanks{Department of Mathematics, The Chinese University of Hong Kong, Shatin, New Territories, Hong Kong SAR, P.R. China (jihohong@cuhk.edu.hk, b.jin@cuhk.edu.hk).}\and Bangti Jin\footnotemark[2]\and Yavar Kian\thanks{Univ Rouen Normandie, CNRS, Normandie Univ, LMRS UMR 6085, F-76000 Rouen, France (\texttt{yavar.kian@univ-rouen.fr})}}

\date{\today}

\begin{document}

\maketitle

\begin{abstract}
In this work we investigate an inverse problem of identifying a spatially variable order in the one-dimensional subdiffusion model from the boundary flux measurement. The model involves a generalized Caputo derivative in time, and arises in the mathematical modeling of anomalous diffusion in heterogeneous media. We prove the unique recovery of a monotone piecewise constant variable order and its range for known and unknown media, respectively. The analysis is based on a delicate asymptotic expansion of the Laplace transform of the data as $p\to0$, which is of independent interest.
\end{abstract}



\section{Introduction}
	
In this work, we investigate an inverse problem for  space-variable order time-fractional diffusion in one spatial dimension.
Fix $L>0$, and let $\alpha:(0,L)\to(0,1)$ be a piecewise constant function. Let $U$ be the solution to the following initial boundary value problem
	\beq\label{eq:withcoeff:ibvp}
	\left\{
	\begin{aligned}
		\rho(x)\p_t^{\alpha(x)}U(t,x) - \p_x( \sigma(x) \p_x U(t,x)) + q(x) U(t,x)&=0,\quad (t,x)\in(0,\infty)\times (0,L),\\
		U(t,0)=g(t),\quad U(t,L)&=0,\quad t\in(0,\infty),\\
		U(0,x)&=0,\quad x\in(0,L),
	\end{aligned}\right.
\eeq
where $\rho\in L^\infty(0,L)$ and $\sigma\in L^\infty(0,L)$ have positive lower bounds, $\sigma$ is Lipschitz continuous, $q\in L^\infty(0,L)$ is nonnegative, and $g$ is the boundary excitation. In the model \eqref{eq:withcoeff:ibvp}, the spatially variable order Caputo fractional derivative $\partial_t^{\alpha(x)} U(t,x)$ in time $t$ is defined by (see, e.g., \cite[p. 92]{KilbasSrivastavaTrujillo:2006} or \cite[p. 41]{Jin:2021})
\begin{equation}
    \partial_t^{\alpha(x)}U(t,x) = \frac{1}{\Gamma(1-\alpha(x))}\int_0^t(t-s)^{\alpha(x)-1}\partial_sU(s,x)\,{\rm d}s,
\end{equation}
where $\Gamma(z)$, for $\Re(z)>0$, denotes Euler's Gamma function. 

Mathematical models of the form \eqref{eq:withcoeff:ibvp} can describe anomalously slow diffusion (subdiffusion) processes. Subdiffusion collectively refers to diffusion phenomena in which the mean squared particle displacement grows sublinearly with the time $t$, and the list of successful practical applications is still growing \cite{MetzlerJeon:2014}. The diffusion process in homogeneous media is frequently modeled by constant-order subdiffusion, in which the mapping $x \mapsto \alpha(x)$ in the model \eqref{eq:withcoeff:ibvp} is taken to be constant over the domain $(0,L)$. This class of mathematical models has been extensively studied both mathematically and numerically. However, in complex media, heterogeneous regions exhibit spatially inhomogeneous variations, and variable-order time-fractional models are deemed more suitable for describing the space-dependent anomalous
diffusion process \cite{SunChen:2009}. Chechkin et al \cite{ChechkinGorenflo:2005} were the first to study a composite system with two regions (i.e., piecewise constant with two pieces) and
found interesting effects involving a nontrivial average drift; see also \cite{KorabelBarkai:2010,Stickler:2011} for related studied in physics. Fedotov and Falconer \cite{Fedotov:2012} exhibited long-time structural instability of variable-order models. 
So far such models have been successfully used in,
e.g., coupled chloride diffusion-binding processes in reinforced concrete \cite{ChenZhang:2013}, liquid infiltration in porous
media \cite{Gerasimov:2010}, spatial heterogeneity in complex anisotropic
medium \cite{Straka:2018}, and nonlinear contact phenomena \cite{Ingman:2000}; 
see also the reviews \cite{SunChangChen:2019,Patnaik:2020}. At a microscopic level, the variable-order subdiffusion corresponds to continuous time random walk (CTRW)-driven
stochastic processes with a space-dependent diffusion coefficient; see the works \cite{ZhangLiLuo:2013,Orsingher:2018,RIcciuti:2017,SavovToaldo:2020} for derivations of variable-order models from CTRW.

In this work, we are concerned with an inverse problem of identifying the fractional order $\alpha$ as a function of the space variable $x$ from the boundary flux data in one space dimension. We will show that under Assumption \ref{ass:main}, the flux measurement $\p_x U(t_k,0)$ for some $\{t_k\}_{k=1}^\infty\subset(0,L)$ uniquely determines a piecewise constant and monotone $\alpha(x)$, and moreover, when the medium properties and boundary excitation are unknown (i.e., $\rho$, $\sigma$, $q$ and $g$ are unknown), then the same measurement can uniquely determine the range of the variable order $\alpha$; see Theorems \ref{thm:Unique:invpb} and \ref{thm:uniqueness:unknowncoeff} in Section \ref{sec:main} for the precise statements. The proofs of these identifiability results are based on a novel first-order asymptotic expansion of the Laplace transforms $\partial_x\widehat{U}(p,0)$ or $\partial_x\widehat{U}(p,L)$ of the solution $U(t,x)$ to \eqref{eq:withcoeff:ibvp} as $p\to 0$, cf. Theorem \ref{thm:asympt:1stord}. The derivation of the asymptotic expansion relies on a delicate analysis of a certain parameterized tridiagonal linear system, and the analysis is of independent interest. The asymptotic analysis utilizes a local solution representation in the Laplace domain between two points of discontinuity of the variable order $\alpha$; see Lemma \ref{lem:laptrans:repr}. This representation involves different tools of the Sturm-Liouville theory  and the fundamental coupling factors $h_j(p)$.
The main challenge lies in the strong coupling between different factors $h_j(p)$ that allows gluing the local representation into a global one, in view of the continuity property \eqref{eq:recursion} of the Laplace transform $\widehat{U}(p,x)$ of the solution  $U(t,x)$ to problem \eqref{eq:withcoeff:ibvp}. To this end, we employ a recursion process involving the terms $r_n$ and $\widetilde r_n$ to represent the key factors $h_j(p)$ in Lemma \ref{lemma:main:coeffrepr}. The analysis of $r_n$ and $\widetilde r_n$ naturally motivates investigating a double sequence $X_l^m$ and its analytic properties, cf. \eqref{eq:def:X} and Lemmas \ref{lemma:X:property} and \ref{lem:X:positive}, which is the main crux of the analysis. We exploit these relations to derive useful integral expressions in the asymptotic expansion of $\partial_x\widehat{U}(p,0)$ as $p\to 0$  in Theorem \ref{thm:asympt:1stord}. These technical details are given in Section \ref{thmsec:asympt:1stord}.

Now we situate this study in the context of inverse problems for subdiffusion. There is a vast literature on inverse problems in the constant-order case (see, e.g., the reviews \cite{JinRundell:2015,LiLiuYamamoto:2019}), but the study in the variable-order case is very limited. Kian, Soccorsi and Yamamoto \cite{Kian:2018:TFD} established the well-posedness of the model using Laplace transform (see also \cite{KS23} for an equivalent definition), and proved the unique recovery of the variable order $\alpha(x)$ from the partial Dirichlet to Neumann map (along with two coefficients). Ikehata and Kian \cite{IkehataKian:2023} proved the unique recovery of the geometry of the variable order from the Neumann data when the problem is equipped with a specially designed Dirichlet data. 
See also \cite{JiangLi:2022,ZhengChengWang:2019} for related results on models involving a time-dependent variable order. By focusing on the 1D case, this work gives two new unique identifiability results for piecewise constant variable orders from one single observation, and complements the works \cite{IkehataKian:2023,Kian:2018:TFD}. 

The rest of the paper is organized as follows. In Section \ref{sec:main}, we state the main results of the work, and provide further discussions. In Section \ref{sec:weak}, we recall the existence and uniqueness of weak solutions to problem \eqref{eq:withcoeff:ibvp}, and derive the solution representation. In Section \ref{thmsec:asympt:1stord}, we present the proof of Theorem \ref{thm:asympt:1stord}, and last in Section \ref{thmsec:Unique:invpb}, we prove the uniqueness results in Theorems \ref{thm:Unique:invpb} and \ref{thm:uniqueness:unknowncoeff}. Throughout, the notation $L^2_\rho(I)$ denotes the space of weighted square integrable functions with a weight $\rho$ over an interval $I$, and $(\cdot,\cdot)_{L_\rho^2(I)}$ denotes the corresponding inner product. 

\section{Main results and discussions}\label{sec:main}

Now we describe the main results of the work. Throughout we make the following assumptions.
\begin{assumption}\label{ass:reg} $\rho\in L^\infty(0,L)$ satisfies $ 0<\rho_0\leq \rho(x) \leq \rho_M <\infty$ for a.e. $x\in (0,L)$, $\sigma\in C^1[0,L]$  satisfies
$0<\sigma_0\leq \sigma(x) \leq \sigma_M <\infty$ for $x\in [0,L]$,
and $q\in L^\infty(0,L)$ is nonnegative.
\end{assumption}

\begin{assumption}\label{ass:main}
The following conditions hold.
\begin{itemize}
\item[{\rm(i)}] The piecewise constant function $\alpha:[0,L]\to (0,1)$ given by $\alpha=\alpha_{0}+\sum_{i=1}^{n} (\alpha_{i}-\alpha_{0})\chi_{(x_{i},x_{i+1})}$ for some partition $x_{0}=0<x_{1}<\dots<x_{n}<x_{n+1}=L$ satisfies
\beq\label{eq:alpha:condition}
	\max_{x\in (0,L)}\alpha(x) < 2\min_{x\in (0,L)}\alpha(x).
\eeq 
\item[{\rm(ii)}] There exist some integer $N\ge2$ and $\{g_k\}_{k=2}^N\subset\mathbb{R}$ with $g_N\ne0$ such that $g$ is given by
\beq\label{eq:g:poly}
    g(t):=\sum_{k=2}^N g_k t^k\quad\mbox{for all }t>0.
\eeq
\end{itemize}
\end{assumption}

We briefly comment on Assumption \ref{ass:main}. Condition (i) is needed also  for the well-posedness of the direct problem \cite[Remark 1]{Kian:2018:TFD} (see also \eqref{est1} and Remark \ref{rmk:well-posed} below) and for deriving the first-order approximation of  $\partial_x\widehat{U}(p,0)$ and $\partial_x\widehat{U}(p,L)$ as $p\to0$ in Theorem \ref{thm:asympt:1stord}. Condition \eqref{eq:alpha:condition} implies that the variable order $\alpha$ (i.e., the sublinearity of the MSD) should not vary greatly at different spatial points, or the heterogeneity of the media is not too strong. For Condition (ii), 
the analysis below indicates that one may choose other analytic function $g$ as long as it can ensure that problem  \eqnref{eq:withcoeff:ibvp} has a unique time-analytic solution with sub-exponential growth in time so that Laplace transform $\widehat{U}(p,x)$ is well defined for all $p>0$. See Theorem \ref{thm:weak} for the unique existence of solutions with sub-exponential growth in time.

First, we state an important first-order asymptotics of the Laplace transform $\partial_x\widehat U(p,0)$ as $p\to 0$ of the solution $U(t,x)$ to problem \eqref{eq:withcoeff:ibvp}. It is  the main technical tool for proving the uniqueness results in Theorems \ref{thm:Unique:invpb} and \ref{thm:uniqueness:unknowncoeff}. The proof is based on a delicate analysis of the solution to a certain parametric tridiagonal linear system, and the details are given in Section \ref{thmsec:asympt:1stord}. 
\begin{thm}\label{thm:asympt:1stord}
Let $\rho$, $\sigma$ and $q$ satisfy Assumption \ref{ass:reg}, and $\alpha$ and $g$ satisfy Assumption \ref{ass:main}.
Then the Laplace transform $\widehat U(p,x)$ of the solution $U(t,x)$ to problem \eqref{eq:withcoeff:ibvp} satisfies 
		\begin{align}
		    \label{eq:asymp:general:duhat}
		\p_x\widehat{U}(p,0)=\widehat{g}(p)\left( C_{0} + \sum_{i=0}^{n} C_{i+1} p^{\alpha_{i}} + O\left(p^{\min\{2\alpha_{i}\,:\,0\le i\le n\}}\right)\right)\quad\mbox{as }p\to0,\\
                \label{eq:asymp:general:duhat:oppside}
		\p_x\widehat{U}(p,L)=\widehat{g}(p)\left( \widetilde{C}_{0} + \sum_{i=0}^{n} \widetilde{C}_{i+1} p^{\alpha_{i}} + O\left(p^{\min\{2\alpha_{i}\,:\,0\le i\le n\}}\right)\right)\quad\mbox{as }p\to0,
		\end{align}
		where $C_{i}$ and $\widetilde{C}_{i}$, $i=0,1,\ldots,n$, are real numbers given by
		\begin{align}
		    \label{eq:C:intform}
		C_{i+1}&=-\|u\|_{L^2_\rho(x_{i},x_{i+1})}^2/\sigma(0) \quad \mbox{and}\quad 
		\widetilde{C}_{i+1}=(u,\overline{u})_{L_\rho^2(x_{i},x_{i+1})}/\sigma(L),
		\end{align}
where	$u$ and $\overline{u}$ are respectively the unique $C^{1,1}$ solutions of the following boundary value problems
		\begin{equation*}
			\left\{\begin{aligned}
				-(\sigma (x)u'(x))' + q(x) u(x)&=0,\quad\mbox{in }(0,L),\\
				u(0)=1,\quad u(L)&=0,
		\end{aligned}\right.
		\end{equation*}
		\begin{equation*}
			\left\{\begin{aligned}
				-(\sigma (x)\overline{u}'(x))' + q(x) \overline{u}(x)&=0,\quad\mbox{in }(0,L),\\
				\overline{u}(0)=0,\quad \overline{u}(L)&=1.
		\end{aligned}\right.
		\end{equation*}
    The functions $u$ and $\overline{u}$ are strictly positive in $(0,L)$ and satisfy
    \beq
    u'(0)=C_0,\quad u'(L)=\widetilde{C}_0\quad\mbox{and}\quad\overline{u}'(0)=-\frac{\sigma(L)}{\sigma(0)}\widetilde{C}_0.
    \eeq
\end{thm}

\begin{remark}
By Tauberian theorem for Laplace transform \cite[Chapter XIII]{Feller:1966}, Theorem \ref{thm:asympt:1stord} yields the large time asymptotics of $\partial_xU(t,0)$ and $\partial_x U(t,L)$. For example, with $s<-1$ and $\widehat{g}(p)=p^{s-1}$, Theorem \ref{thm:asympt:1stord} implies 
\begin{align*}
	\p_xU(t,0)&=\frac{C_{0}t^{-s}}{\Gamma(1-s)}  + \sum_{i=0}^{n} \frac{C_{i+1} t^{-s-\alpha_{i}}}{\Gamma(1-s-\alpha_i)} + O\left(t^{-s-\min\{2\alpha_{i}\,:\,0\le i\le n\}}\right),\quad\mbox{as }t\to\infty,\\
		\p_xU(t,L)&=\frac{\widetilde{C}_0t^{-s}}{\Gamma(1-s)}+ \sum_{i=0}^{n} \frac{\widetilde{C}_{i+1}t^{-s-\alpha_{i}}}{\Gamma(1-s-\alpha_i)} + O\left(t^{-s-\min\{2\alpha_{i}\,:\,0\le i\le n\}}\right),\quad\mbox{as }t\to\infty.
\end{align*}
\end{remark}

Now we state the setting for the inverse problem. For $j=1,2$, let $L^j\in(0,\infty)$, let $q^j\in L^\infty(0,L^j)$ be nonnegative, let $\rho^j\in L^\infty(0,L^j)$ satisfy $0<\rho_0\leq \rho^j(x)\leq \rho_M$ for $x\in (0,L^j)$, {\color{blue}let} $g$ satisfy condition \eqnref{eq:g:poly}, let $\alpha^j$ be piecewise constant and satisfy condition \eqnref{eq:alpha:condition} and let
$U^j$ be the solution of problem \eqnref{eq:withcoeff:ibvp} corresponding to  $(L,q,\rho,\alpha)=(L^j,q^j,\rho^j,\alpha^j)$.
Let $x_{0}^j=0<x_{1}^j<\dots<x_{n^j}^j<x_{n^j+1}^j=L^j$ and
\beq\label{eq:alpha:pc}
\alpha^j=\alpha_{0}^j+\sum_{i=1}^{n^j} \big(\alpha_{i}^j-\alpha_{0}^j\big)\chi_{(x_{i}^j,x_{i+1}^j)},\quad\mbox{with }j=1,2.
\eeq
The next result gives the main result for the inverse problem, i.e., unique identifiability of the order $\alpha$ from the boundary flux data $\partial_xU(t_k,0)$  when the excitation $g$ is given at either $x=0$ or $x=L$. The proof of Theorem \ref{thm:Unique:invpb} is given in Section \ref{thmsec:Unique:invpb}. The key of the proof is the following observation from Theorem \ref{thm:asympt:1stord}: $L$ and $u$ [$\overline{u}$] are uniquely determined by $\sigma$, $q$ and $C_0$ [$\widetilde{C}_0$].
\begin{thm}[Unique identifiability]\label{thm:Unique:invpb}
For $j=1,2$, let $L^j\in(0,\infty)$, $\alpha^j$ be piecewise constant and monotone increasing {\rm(}respectively decreasing{\rm)} and satisfy condition \eqref{eq:alpha:condition}, and let $g$ satisfy \eqref{eq:g:poly}.
Let $\rho$, $\sigma$ and $q$ satisfy Assumption \ref{ass:reg}.
Let $U^j$ be the solution of problem \eqref{eq:withcoeff:ibvp} with $(L,\alpha)=(L^j,\alpha^j)$ for $j=1,2$ and one of the following boundary conditions:
\begin{itemize}
    \item[{\rm(i)}] for $j=1,2$, $U^j(t,0)=g(t)$ and $U^j(t,L^j)=0$ for $t\in(0,\infty)$; 
    \item[{\rm(ii)}] for $j=1,2$, $U^j(t,0)=0$ and $U^j(t,L^j)=g(t)$ for $t\in(0,\infty)$.
\end{itemize}
If $\p_x U^1(t_k,0)=\p_x U^2(t_k,0)$ for some sequence $\{t_k\}_{k=1}^\infty$ of distinct numbers converging to a positive number, then we have $L^1=L^2$ and $\alpha^1=\alpha^2$.
\end{thm}
  
Then we consider the case with unknown excitation $g$ and coefficients $\rho$, $\sigma$ and $q$, using Theorem \ref{thm:asympt:1stord}
and the universal properties $C_{i+1}^j<0$ and $\widetilde{C}_{i+1}^j>0$ for all $i\ge0$ and $j=1,2$. Theorem \ref{thm:uniqueness:unknowncoeff} indicates that even for an unknown excitation $g$ and coefficients $\rho$, $\sigma$ and $q$, the range of the piecewise-constant variable order $\alpha$ is uniquely determined. It is worth emphasizing that the result does not require the variable order $\alpha$ to be monotone. Theorem \ref{thm:uniqueness:unknowncoeff} generalizes existing results on determining a constant order, multiple orders and distributed order in an unknown medium \cite{JinKian:2022siap,JinKian:2021procA,JinKian:2023cms}.

\begin{thm}[Uniqueness of range for unknown coefficients] \label{thm:uniqueness:unknowncoeff}
For $j=1,2$, let $L^j\in(0,\infty)$, let $(\rho,\sigma,q)=(\rho^j,\sigma^j,q^j)$ satisfy Assumption \ref{ass:reg} with $L=L^j$, let $\alpha^j:(0,L^j)\to(0,1)$ be piecewise constant and satisfy condition \eqref{eq:alpha:condition}, and let $g^j$ satisfy condition \eqref{eq:g:poly} for some $N^j$ and $\{g_{k}^j\}_{k=2}^{N^j}$. Let $U^j$ with $j=1,2$ be the solution of
\beq\label{eq:withcoeff:ibvp:twocases}
\left\{
\begin{aligned}	\rho^j(x)\p_t^{\alpha^j(x)}U^j(t,x) - \p_x( \sigma^j(x) \p_x U^j(t,x)) + q^j(x) U^j(t,x)&=0,\quad (t,x)\in(0,\infty)\times (0,L^j),\\
U^j(t,a^j)=g_j(t),\quad U^j(t,b^j)&=0,\quad t\in(0,\infty),\\
U^j(0,x)&=0,\quad x\in(0,L^j),
	\end{aligned}\right.
\eeq
where $(a^j,b^j)=(0,L^j)$ or $(a^j,b^j)=(L^j,0)$. If $\p_x U^1(t_k,0)=\p_x U^2(t_k,0)$ for some sequence $\{t_k\}_{k=1}^\infty$ of distinct numbers converging to a positive number, then we have
\beq\label{eq:alpha:identicalset}
\{\alpha^1(x)\,:\,x\in(0,L^1)\}=\{\alpha^2(x)\,:\,x\in(0,L^2)\}.
\eeq
\end{thm}	

\section{Preliminaries: Existence, uniqueness and representation}\label{sec:weak}
\subsection{Existence and uniqueness of solutions}

First we discuss the existence and uniqueness of strong solutions to problem \eqref{eq:withcoeff:ibvp} in the following sense.
\begin{definition}\label{d1}
A function $U\in C^1([0,\infty);L^2(0,L))\cap C([0,\infty);H^2(0,L))$ is said to be a strong solution of problem \eqref{eq:withcoeff:ibvp} if the following two conditions are fulfilled: {\rm(a)} for all $\psi\in L^2(0,L)$,  we have
$$\int_0^L(\rho(x)\p_t^{\alpha(x)}U(t,x) - \p_x( \sigma(x) \p_x U(t,x)) + q(x) U(t,x))\psi(x)\,{\rm d} x=0,\quad \forall t\in(0,\infty);$$
and {\rm(b)} $U(t,0)= g(t)$, $U(t,L)= 0$,   for all  $t\in(0,\infty)$, 
and  $U(0,\cdot)\equiv0$.
\end{definition}

\begin{thm}\label{thm:weak}
Let $\rho$, $\sigma$ and $q$ satisfy Assumption \ref{ass:reg}, and $\alpha$ and $g$ satisfy Assumption \ref{ass:main}.
Then there is a unique strong solution $U$ of problem \eqref{eq:withcoeff:ibvp} with sub-exponential growth in time. Moreover, the map $t\mapsto U(t,\cdot)$ is analytic with respect to $t\in(0,\infty)$ as a map valued in $H^2(0,L)$.
\end{thm}
\begin{proof}
The proof follows closely \cite[Proposition 3.1]{Kian:2018:TFD}. Let  $\underline{\alpha}=\min_{x\in(0,L)}\alpha(x)$ and $\overline{\alpha}=\max_{x\in(0,L)}\alpha(x)$, and let $A$ be the Dirichlet realization of the operator 
$- \p_x( \sigma(x) \p_x  ) + q(x)$
acting on $L^2(0,L)$ with its domain $D(A)= H^1_0(0,L)\cap H^2(0,L)$. By \cite[Proposition 2.1]{Kian:2018:TFD}, for all $p\in\mathbb C\setminus\mathbb R_-$, $A+\rho(x)p^{\alpha(x)}$ is boundedly invertible in $L^2(0,L)$, $(A+\rho(x)p^{\alpha(x)})^{-1}$ maps $L^2(0,L)$ to $H^2(0,L)\cap H^1_0(0,L)$ and for all $\psi\in(0,\pi)$, 
\begin{equation}\label{est1}\|(A+\rho(x)(re^{i\phi})^{\alpha(x)})^{-1}\|_{\mathcal B(L^2(0,L))}\leq C\max(r^{\underline{\alpha}-2\overline{\alpha}},r^{\overline{\alpha}-2\underline{\alpha}}),\quad r>0,\ \phi\in(-\psi,\psi),
\end{equation}
where $C>0$ is independent of $r$, and $\|\cdot\|_{\mathcal B(L^2(0,L))}$  denotes the operator norm from $L^2(0,L)$ to itself. For all $p\in\mathbb C\setminus\mathbb R_-$, consider the boundary value problem
\beq\label{eqq1}
    \left\{\begin{aligned}
	-\p_x(\sigma(x)\p_x\widehat V(p,x))+ (\rho(x) p^{\alpha(x)} + q(x)) \widehat V(p,x)&=0,\quad\mbox{in } (0,L),\\
	\widehat V(p,0)=\widehat{g}(p), \quad \widehat V(p,L)&=0.\end{aligned}\right.
\eeq
We can decompose $\widehat V(p,\cdot)$ into $\widehat V(p,x)=\widehat{g}(p)h(x)+\widehat W(p,x)$ for $x\in(0,L),$
with
$\widehat W(p,\cdot)=-\widehat{g}(p)(A+\rho(x)p^{\alpha(x)})^{-1}\rho(x)p^{\alpha(x)}h$ and $h\in H^2(0,L)$ satisfying
$$
\left\{\begin{aligned}
			-\p_x(\sigma(x)\p_xh(x))  + q(x)h(x)&=0,\quad\mbox{in } (0,L),\\
			h(0)=1, \quad h(L)&=0.\end{aligned}\right.$$
Fix
$\theta\in(\frac{\pi}{2},\pi)$, $\delta>0$ and define a contour  
$\gamma(\delta,\theta):=\gamma_-(\delta,\theta)\cup\gamma_0(\delta,\theta)\cup\gamma_+(\delta,\theta)\subset\mathbb{C}$, oriented  counterclockwise with
$\gamma_0(\delta,\theta):=\{\delta\, e^{i\beta}:\ \beta\in[-\theta,\theta]\}$, $\gamma_\pm(\delta,\theta)
:=\{s\,e^{\pm i\theta}: s\in[\delta,\infty)\}$.
Since $\widehat{g}(p)=\sum_{k=2}^N g_kk!p^{-k-1}$
and by applying \eqref{est1}, similar to \cite[Proposition 3.1]{Kian:2018:TFD}, we can define the map
$$U(t,\cdot)=g(t)h+\frac{1}{2i\pi}\int_{\gamma(\delta,\theta)}e^{t p} \widehat W(p,\cdot)\,{\rm d} p,\quad t\in[0,\infty).$$
The Laplace transform  $\widehat{U}(p,\cdot)\in H^2(0,L)$ of $U(t,\cdot)$ is defined for all $p>0$ and satisfies $\widehat{U}(p,\cdot)=\widehat V(p,\cdot)$. Thus the desired assertion holds if $U\in C^1([0,\infty);L^2(0,L))\cap C([0,\infty);H^2(0,L))$.
It suffices to prove $$W(t,\cdot) = \frac{1}{2i\pi}\int_{\gamma(\delta,\theta)}e^{t p} \widehat W(p,\cdot)\,{\rm d} p \in  C^1([0,\infty);L^2(0,L))\cap C([0,\infty);H^2(0,L)\cap H^1_0(0,L)).$$ 
The resolvent estimate \eqref{est1} implies that for all $r\geq1$, there exists $C>0$ independent of $r$ such that
$$\begin{aligned} &\|\widehat W(re^{\pm i\theta},\cdot)\|_{D(A)}\leq C\|A\widehat W(re^{\pm i\theta},\cdot)\|_{L^2(0,L)}
\leq Cr^{3\overline{\alpha}-2\underline{\alpha}-3}.\end{aligned}$$
Now condition \eqref{eq:alpha:condition} implies $3\overline{\alpha}-2\underline{\alpha}-3<2\overline{\alpha}-3<-1$. Hence, we obtain $r\mapsto \widehat W(re^{\pm i\theta},\cdot)\in L^1(\delta,\infty;H^2(0,L)\cap H^1_0(0,L))$ and deduce
\begin{equation}\label{rep1} W(t,\cdot)\in C([0,\infty);H^2(0,L)\cap H^1_0(0,L)).\end{equation}
Similarly, for all $r\geq1$, we find
$$ \|re^{\pm i\theta}\widehat W(re^{\pm i\theta},\cdot)\|_{L^2(0,L)}\leq C r^{2\overline{\alpha}-2\underline{\alpha}-2},$$
and, since  $2\overline{\alpha}-2\underline{\alpha}-2<\overline{\alpha}-2<-1$, we have $r\mapsto re^{\pm i\theta}\widehat W(re^{\pm i\theta},\cdot)\in L^1(\delta,\infty;L^2(0,L))$ which in turn implies 
$W(t,\cdot)\in C^1([0,\infty);L^2(0,L)).$
This proves that $U\in C^1([0,\infty);L^2(0,L))\cap C([0,\infty);H^2(0,L))$ with
\begin{equation}\label{rep2}U(t,\cdot)=g(t)h+W(t,\cdot),\quad t\in[0,\infty).
\end{equation}
Note that
$U(0,\cdot)=\frac{1}{2i\pi}\int_{\gamma(\delta,\theta)} \widehat W(p,\cdot){\rm d} p.$
By \cite[Proposition 3.1]{Kian:2018:TFD}, the map $p\mapsto (A+\rho(x)p^{\alpha(x)})^{-1}$ is holomrphic on $\mathbb C\setminus(-\infty,0]$ as an operator valued in $\mathcal B(L^2(0,L))$, which implies that the map $p\mapsto \widehat W(p,\cdot)$ is also holomorphic on $\mathbb C\setminus(-\infty,0]$ as a map valued in $L^2(0,L)$.
Therefore, similar to \cite[p. 16]{Ki23}, by fixing $R>\delta$  and applying the residue theorem, we deduce
\begin{align*}
\frac{1}{2i\pi}\int_{\gamma(\delta,\theta)} \widehat W(p,\cdot)\,{\rm d} p=\lim_{R\to\infty}\frac{1}{2i\pi}\int_{ \gamma_0(R,\theta)} \widehat W(p,\cdot)\,{\rm d}p.
\end{align*}
Meanwhile, in view of the resolvent estimate \eqref{est1}, we have
\begin{equation*}
\|Re^{i\beta} \widehat W(Re^{i\beta},\cdot)\|_{L^2(0,L)}\leq CR^{2\overline{\alpha}-2\underline{\alpha}-2},\quad R>1,\ \beta\in[-\theta,\theta],
\end{equation*}
which implies
\begin{equation*}
\frac{1}{2i\pi}\int_{\gamma(\delta,\theta)} \widehat W(p,\cdot){\rm d} p = \lim_{R\to\infty}\frac{1}{2i\pi}\int_{\gamma_0(R,\theta)} \widehat W(p,\cdot){\rm d} p\equiv 0.
\end{equation*}
Thus we have $U(0,\cdot)\equiv 0$, and in view of \eqref{rep1}-\eqref{rep2}, we have
$$U(t,0)=g(t)h(0)=g(t),\quad U(t,1)=g(t)h(1)=0,\quad t\in[0,\infty).$$
This proves that $U$ fulfills condition (b) of Definition \ref{d1}. Next, fix
$$G(t,x)=\rho(x)\p_t^{\alpha(x)}U(t,x) - \p_x( \sigma(x) \p_x U(t,x)) + q(x) U(t,x),\quad (t,x)\in[0,\infty)\times(0,L)$$
and note that $G\in C([0,\infty);L^2(0,L))$ admits sub-expontial growth in time. Since $U(0,\cdot)\equiv 0$,   for all $p>0$, Laplace transform in time $\widehat{G}(p,\cdot)$ of $G$ is given by
$$\widehat{G}(p,x)=-\p_x(\sigma(x)\p_x\widehat{U}(p,x))+ \big(\rho(x) p^{\alpha(x)} + q(x)\big) \widehat{U}(p,x),\quad x\in(0,L).$$
Since $\widehat{U}(p,\cdot)=\widehat V(p,\cdot)$ solves \eqref{eqq1}, we deduce that $\widehat{G}(p,\cdot)\equiv0$, for $p>0$, and by the injectivity of Laplace transform, it follows that $G\equiv0$. Thus, condition (a) of Definition \ref{d1} is also fulfilled and $U$ is a strong solution of \eqref{eq:withcoeff:ibvp}. The uniqueness follows from that of Laplace transform which is necessarily the unique solution of  problem \eqref{eqq1} (see, e.g., the proof of \cite[Theorem 1.3]{Ki23}). The analyticity of the map $t\mapsto U(t,\cdot)$ with respect to $t\in(0,\infty)$ as a map valued in $H^2(0,L)$ can be deduced by combining the above argumentation with that of \cite[Lemma 3.2]{Kian:2018:TFD}.
\end{proof}

\begin{remark}
Compared with existing works \cite{Kian:2018:TFD,Ki23}, there are two new aspects of Theorem \ref{thm:weak}: it gives the $C^1$ regularity in time and verifies the initial condition in a strong sense.    
\end{remark}
\begin{remark}\label{rmk:well-posed} Condition \eqref{eq:alpha:condition} ensures that the solution $u$ to problem \eqref{eq:withcoeff:ibvp} can be represented via Duhamel's principle \cite[Remark 1]{Kian:2018:TFD}. In the analysis, problem \eqref{eq:withcoeff:ibvp} is stated for $t\in(0,\infty)$. This is due to the use of Laplace transform to establish its well-posedness, which, to the best of our knowledge, is the only known analysis technique for proving the well-posedness of space variable-order time-fractional models \cite{Kian:2018:TFD}. For a problem defined on a finite time interval $(0,T)$, the solution $u$ can be obtained by first suitably extending the boundary data from $(0,T)$ to the infinite time interval $(0,\infty)$ and then restricting the corresponding solution $u$ on $(0,\infty)$ to $(0,T)$. 
\end{remark}

Below we derive an expression for the Laplace transform $\widehat{U}=\widehat{U}(p,x)$ of the solution $U=U(t,x)$ to problem \eqref{eq:withcoeff:ibvp}.
		Using the sub-exponential growth of $U(t,x)$ in time $t$, cf. Theorem \ref{thm:weak}, the Laplace transform $\widehat{U}(p,\cdot)\in H^2(0,L)$ of $U(t,\cdot)$ is defined for all $p>0$ and satisfies
		\beq
		\left\{\begin{aligned}
			-\p_x(\sigma(x)\p_x\widehat{U}(p,x))+ \big(\rho(x) p^{\alpha(x)} + q(x)\big) \widehat{U}(p,x)&=0,\quad\mbox{in } (0,L),\\
			\widehat{U}(p,0)=\widehat{g}(p), \quad \widehat{U}(p,L)&=0.\end{aligned}\right.
		\eeq

\subsection{Eigenfunction expansion}
We derive an explicit expression of $\widehat{U}(p,x)$ using the Strum-Liouville theory \cite{Poschel:1987:IST}. Fix $j=0,1,\dots,n$.
We define the elliptic operator $A_j$ on the Hilbert space $L_\rho^2(x_j,x_{j+1})$ by
$$A_j [f](x) := \rho(x)^{-1}(-\p_x(\sigma(x)\p_xf(x)) + q(x)f(x))\quad\mbox{for }x\in(x_j,x_{j+1})\mbox{ and }f\in D(A_j),$$
where its domain $D(A_j):=H_0^1(x_j,x_{j+1})\cap H^2(x_j,x_{j+1})$ is dense in the space $L_\rho^2(x_j,x_{j+1})$.
Let $\{(\lambda_{j,k},\varphi_{j,k})\}_{k=1}^\infty$ be the Dirichlet eigenpairs of $A_j$, with the eigenvalues $\{\lambda_{j,k}\}_{k=1}^\infty$ ordered nondecreasingly, normalized by $\|\varphi_{j,k}\|_{L_\rho^2(x_j,x_{j+1})}=1$.
    From \cite[Theorem 2, Chapter 2]{Poschel:1987:IST}, all the eigenvalues are positive and simple.
    Since $A_j$ is self-adjoint, $\{\varphi_{j,k}\}_{k=1}^\infty$ forms an orthonormal basis of the space $L_\rho^2(x_j,x_{j+1})$.
	Let $v_j$ and $w_j$ be the fundamental solutions of the elliptic operator $-\p_x(\sigma(x)\p_x) + q(x)$, defined respectively by
		\begin{align*}
			\left\{\begin{aligned}
				-(\sigma(x)v_j'(x))'+ q(x) v_j(x)&=0,\quad\mbox{in } (x_j,x_{j+1}),\\
				v_j(x_j)=1, \quad v_j(x_{j+1})&=0,\end{aligned}\right.\\
			\left\{\begin{aligned}
				-(\sigma(x)w_j'(x))' + q(x) w_j(x)&=0,\quad\mbox{in } (x_j,x_{j+1}),\\
				w_j(x_j)=0, \quad w_j(x_{j+1})&=1.\end{aligned}\right.
		\end{align*}
		Using integration by parts, direct computation gives, for $j=0,1,\dots,n$ and $k\in\mathbb{N}$,
		\beq\label{eq:vw:phicoeff}
		\begin{aligned}
			&(v_j,\varphi_{j,k})_{L_\rho^2(x_{j},x_{j+1})} = \frac{{\sigma(x_j)}\varphi_{j,k}'(x_{j})}{\lambda_{j,k}}\quad\mbox{and}\quad (w_j,\varphi_{j,k})_{L_\rho^2(x_{j},{x_{j+1}})} = -\frac{{\sigma(x_{j+1})}\varphi_{j,k}'(x_{j+1})}{\lambda_{j,k}}.
		\end{aligned}
		\eeq
Then we have the following solution representation of the Laplace transform $\widehat U(p,x)$. This representation plays an important role in deriving the first-order approximation in Section \ref{thmsec:asympt:1stord}. However, it involves the strong coupling between the factors $h_j(p)$, which greatly complicates the subsequent analysis.

\begin{lemma}\label{lem:laptrans:repr}
For $j=0,1,\dots,n$, $x\in(x_j,x_{j+1})$ and $p>0$, we have
\begin{align}
	\widehat{U}(p,x)&=h_j(p)v_j(x) + h_{j+1}(p)w_j(x) \nonumber\\
    &\quad + \sum_{k=1}^\infty \frac{1}{\lambda_{j,k}}\frac{p^{\alpha_j}}{p^{\alpha_j}+\lambda_{j,k}}\left[h_{j+1}(p)\sigma(x_{j+1})\varphi_{j,k}'(x_{j+1})-h_j(p)\sigma(x_j)\varphi_{j,k}'(x_j)\right]\varphi_{j,k}(x),\label{eq:Uhat:eigftnexp}
\end{align}
where $h_j(p)$ is defined by $h_j(p):=\lim_{x\to x_j} \widehat{U}(p,x)$, $j=1,\ldots,n$, with  $h_0(p)=\widehat{g}(p)$ and $h_{n+1}(p)=0$.
Moreover, the coefficient functions $h_j$ satisfy	\begin{equation}\label{eq:recursion}
	\lim_{s\to 0^+}\big(\p_x\widehat{U}(p,x_j+s)-\p_x\widehat{U}(p,x_j-s)\big)=0,\quad j=1,\dots,n.
\end{equation}
\end{lemma}
\begin{proof}
Using the sub-exponential growth of $U(t,x)$ in time $t$ from Theorem \ref{thm:weak}, the Laplace transform $\widehat{U}(p,\cdot)\in H^2(0,L)$ of $U(t,\cdot)$ is defined for all $p>0$, and thus \eqref{eq:recursion} holds. For any $p>0$, we set 
\begin{align*}
	\widehat W_j(p,x)&:=\widehat{U}(p,\cdot) - h_j(p) v_j(x) - h_{j+1}(p) w_j(x),\quad x\in(x_j,x_{j+1}).	
\end{align*}
For any $p>0$, the function $\widehat{W}_j$ satisfies 
\begin{equation}\label{eq:afterLaptrans:1D}
   \left\{\begin{aligned}
	    \left(-\p_x(\sigma\p_x) + \rho p^{\alpha_j} + q\right)\widehat W_j(p,x)&=-\rho p^{\alpha_j}\left(h_j(p)v_j+h_{j+1}(p)w_j\right),\quad \mbox{in }(x_j,x_{j+1}),\\
	\widehat W_j(p,x_j)=\widehat W(p,x_{j+1})&=0,
   \end{aligned}\right.
\end{equation}
and admits a basis expansion with respect to $\{\varphi_{j,k}\}_{j=1}^\infty$ in the  $L_\rho^2(x_j,x_{j+1})$ inner product:
$$\widehat W_j(p,x) = \sum_{k=1}^\infty B_{j,k}(p) \varphi_{j,k}(x),\quad\mbox{with}\quad B_{j,k}(p):=(\widehat W_j(p,\cdot),\varphi_{j,k})_{L^2_\rho(x_j,x_{j+1})}.$$
By integration by parts and using the governing equation for $\widehat{W}_j(p,x)$ in \eqref{eq:afterLaptrans:1D}, we obtain
\begin{align*}
	B_{j,k}(p)&=\int_{x_j}^{x_{j+1}}\widehat W_j(p,x)\frac{1}{\lambda_{j,k}}\left(-(\sigma\varphi_{j,k}')'(x) + q\varphi_{j,k}(x)\right)\,{\rm d}x\\
    &=\frac{1}{\lambda_{j,k}}\int_{x_j}^{x_{j+1}}\left(-(\sigma \widehat W_j')'(p,x) + q\widehat W_j(p,x)\right)\varphi_{j,k}(x)\,{\rm d}x\\
	&=-\frac{p^{\alpha_j}}{\lambda_{j,k}}\left[B_{j,k}(p) + \int_{x_j}^{x_{j+1}}\left(h_j(p)v_j(x)+h_{j+1}(p)w_j(x)\right)\varphi_{j,k}(x)\rho(x)\,{\rm d}x\right].
\end{align*}
Using the identities in \eqnref{eq:vw:phicoeff}, we arrive at the desired representation \eqnref{eq:Uhat:eigftnexp}.
\end{proof}
		
\subsection{Determination of the factors $h_1$ and $h_n$}
The coefficients $h_j$, $j=1,\ldots,n$, of the expression \eqref{eq:Uhat:eigftnexp} are determined by the system \eqref{eq:recursion}. For the following analysis, we define three auxiliary functions $E_j(p), F_j(p)$ and $G_j(p)$ respectively by
\begin{align*}
	E_j(p)&:=\sum_{k=1}^\infty \frac{\big(\sigma(x_j)\varphi_{j,k}'(x_j)\big)^2}{\lambda_{j,k}}\frac{p^{\alpha_j}}{p^{\alpha_j}+\lambda_{j,k}},\\
	F_j(p)&:=\sum_{k=1}^\infty \frac{\sigma(x_{j+1})\sigma(x_j)\varphi_{j,k}'(x_{j+1})\varphi_{j,k}'(x_{j})}{\lambda_{j,k}}\frac{p^{\alpha_j}}{p^{\alpha_j}+\lambda_{j,k}},\\
 	G_j(p)&:=\sum_{k=1}^\infty \frac{\big(\sigma(x_{j+1})\varphi_{j,k}'(x_{j+1})\big)^2}{\lambda_{j,k}}\frac{p^{\alpha_j}}{p^{\alpha_j}+\lambda_{j,k}}.
\end{align*}
These functions will be used to simplify the identities \eqref{eq:recursion}; see the identities \eqref{eq:recursion:explicit} below.
\begin{lemma}
The three functions $E_j(p)$, $F_j(p)$ and $G_j(p)$ are absolutely convergent for $p>0$.
\end{lemma}
\begin{proof}
For the proof, we set
\begin{align*}
E_{j}^*&:=\sum_{k=1}^\infty\frac{\big(\sigma(x_j)\varphi_{j,k}'(x_j)\big)^2}{\lambda_{j,k}^2},\quad F_{j}^*:=\sum_{k=1}^\infty\frac{\sigma(x_{j+1})\sigma(x_j)\varphi_{j,k}'(x_j)\varphi_{j,k}'(x_{j+1})}{\lambda_{j,k}^2},\\ G_{j}^*&:=\sum_{k=1}^\infty\frac{\big(\sigma(x_{j+1})\varphi_{j,k}'(x_{j+1})\big)^2}{\lambda_{j,k}^2}.
	\end{align*}
Since $\{\varphi_{j,k}\}_{k=1}^\infty$ is an orthonormal basis of $L_\rho^2(x_{j},x_{j+1})$ and \eqnref{eq:vw:phicoeff} holds, by Parseval's identity, we have 
\beq\label{eq:E:Parseval}
     \begin{aligned}
	E_{j}^*&=\|v_j\|_{L^2_\rho(x_{j},{x_{j+1}})}^2  ,\quad
	F_{j}^*=-(v_j,w_j)_{L^2_\rho(x_{j},x_{j+1})} ,\quad
	G_{j}^*=\|w_j\|_{L^2_\rho(x_{j},x_{j+1})} ^2.
\end{aligned}
\eeq
Meanwhile, from the inequality
$\frac{p^{\alpha_j}}{p^{\alpha_j}+\lambda_{j,k}}\le\frac{p^{\alpha_j}}{\lambda_{j,k}}$,
we have, for all $j=0,1,\dots,n$ and $p>0$,
\begin{equation}\label{eq:EG:bounds}
	0\le E_j(p)\le p^{\alpha_j}E_j^*\quad\mbox{and}\quad 0\le G_j(p)\le p^{\alpha_j}G_j^*.
\end{equation}
Thus, we have $E_j(p)<\infty$ and $G_j(p)<\infty$ for all $p>0$. The series $F_j(p)$ is also absolutely convergent because the absolute series is bounded by $\frac{1}{2}\left(E_j(p)+G_j(p)\right)$.    \end{proof}

Now we use the functions $E_j$, $F_j$ and $G_j$ to study the system \eqref{eq:recursion}. It follows from  \eqref{eq:Uhat:eigftnexp} that
\begin{align*}	&\lim_{s\to0^+}\sigma(x_j+s)\partial_x\widehat{U}(p,x_j+s)
=h_j(p)\sigma(x_j)v_j'(x_j) + h_{j+1}(p)\sigma(x_j)w_j'(x_j)\\
	&\quad + \sum_{k=1}^\infty \frac{1}{\lambda_{j,k}}\frac{p^{\alpha_j}}{p^{\alpha_j}+\lambda_{j,k}}\left[h_{j+1}(p)\sigma(x_{j+1})\varphi_{j,k}'(x_{j+1})-h_j(p)\sigma(x_j)\varphi_{j,k}'(x_j)\right]\sigma(x_j)\varphi_{j,k}'(x_j)\\
&=h_j(p)\sigma(x_j)v_j'(x_j) + h_{j+1}(p)\sigma(x_j)w_j'(x_j)  + h_{j+1}(p)F_j(p)-h_j(p)E_j(p).
\end{align*}
Similarly,
\begin{align*}
&\lim_{s\to0^+}\sigma(x_{j}-s)\partial_x\widehat{U}(p,x_{j}-s)\\
=&h_{j-1}(p)\sigma(x_{j})v_{j-1}'(x_{j}) + h_{j}(p)\sigma(x_{j})w_{j-1}'(x_{j}) +h_{j}(p)G_{j-1}(p)-h_{j-1}(p)F_{j-1}(p).
\end{align*}
Thus, the identities \eqref{eq:recursion} are equivalent to the following tridiagonal linear system
\begin{align}
	&h_j(p)\left(E_j(p) + G_{j-1}(p) - \sigma(x_j)v_j'(x_j) + \sigma(x_j)w_{j-1}'(x_j)  \right) \nonumber\\
	=& h_{j-1}(p)\left(F_{j-1}(p) - \sigma(x_j)v_{j-1}'(x_j) \right) + h_{j+1}(p)\left(F_j(p) + \sigma(x_j)w_j'(x_j)\right),\quad j=1,\dots,n.\label{eq:recursion:explicit}
\end{align}
The goal is to derive the asymptotics of the solution components $h_1$ and $h_n$ of the linear system \eqnref{eq:recursion:explicit}. This will be conducted in Section \ref{thmsec:asympt:1stord}. Note that the factors $h_0$ and $h_{n+1}$ are given by the boundary conditions of $U(t,x)$. Then using the representation \eqref{eq:Uhat:eigftnexp}, one may obtain the boundary flux data. 

It follows from \eqnref{eq:EG:bounds} that all three functions $E_j(p)$, $F_j(p)$ and $G_j(p)$ tend to zero as $p\to0$. Since $\sigma(x_j)v_{j-1}'(x_j)\ne0$ holds (see, e.g., Lemma \ref{lem:v:decreasing} below), there exists a $p_0>0$ such that
\begin{equation*}
F_{j-1}(p) - \sigma(x_j)v_{j-1}'(x_j)\ne0,\quad\forall p\in(0,p_0)\mbox{ and }j=1,\dots,n.
\end{equation*}
For notational simplicity, we set $c_m=c_m(p)$ and $d_m=d_m(p)$ for $m=1,\dots,n$ as
\beq\label{def:cmdm}
\begin{aligned}
	c_m(p)&:=\frac{E_m(p)+G_{m-1}(p)-\sigma(x_m)v_m'(x_m)+\sigma(x_m)w_{m-1}'(x_{m})}{F_{m-1}(p)-\sigma(x_m)v_{m-1}'(x_m)},\\ d_m(p)&:=-\frac{F_m(p)+\sigma(x_m)w_m'(x_{m})}{F_{m-1}(p)-\sigma(x_m)v_{m-1}'(x_m)},
		\end{aligned}
\eeq
for all $p\in(0,p_0)$. Then, for each $p\in(0,p_0)$, the tridiagonal linear system \eqnref{eq:recursion:explicit} can be written into
\begin{equation}\label{eq:recursion:hcd}
		h_{m-1}(p)=c_m(p)h_m(p) + d_m(p) h_{m+1}(p),\quad m=1,\dots,n.
\end{equation}
The coupling between different components poses the main challenge in deriving the asymptotic expansion of $h_1(p)$ and $h_n(p)$, which are directly connected with the boundary measurement. To this end, we define $r_m=r_m(p)$ and $\widetilde{r}_m=\widetilde{r}_m(p)$ for $p\in(0,p_0)$ and $m=1,\dots,n$ recursively by
		\beq\label{eq:recursion:rns}
		\begin{aligned}
			&r_0=\widetilde{r}_1=1,\quad s_0=\widetilde{s}_1=0,\\
			&r_m=c_m r_{m-1} + s_{m-1},\quad s_m=d_m r_{m-1},\quad\mbox{for }m=1,2,\dots,n,\\
			&\widetilde{r}_m=c_m \widetilde{r}_{m-1} + \widetilde{s}_{m-1},\quad \widetilde{s}_m=d_m \widetilde{r}_{m-1},\quad\mbox{for }m=2,\dots,n.
		\end{aligned}
		\eeq
				
The next result gives key identities connecting $h_0(p)$ and $h_1(p)$ with $r_n(p)$ and $\widetilde{r}_n(p)$.
\begin{lemma}\label{lemma:main:coeffrepr}
The following two identities hold
\begin{align*} r_n(p) h_n(p) = h_0(p)\quad \mbox{and}\quad  \widetilde{r}_n(p) h_n(p) = h_1(p),\quad \forall p\in(0,p_0).
\end{align*}
\end{lemma}
\begin{proof}
Fix $p\in(0,p_0)$. Using mathematical induction on $k$, we will prove
\begin{equation}\label{eq:rnhn:obj}
    r_nh_n = r_{n-k}h_{n-k} + s_{n-k}h_{n-k+1},\quad k=1,\dots,n.
\end{equation}
By the recursion \eqnref{eq:recursion:rns} with $m=n$, the identity \eqnref{eq:recursion:hcd} with $m=n$ and the condition $h_{n+1}=0$, we have
\begin{equation*}
r_nh_n = (c_n r_{n-1} + s_{n-1})h_n = r_{n-1}h_{n-1} + s_{n-1}h_{n}.
\end{equation*}
That is, the claim \eqref{eq:rnhn:obj} holds for $k=1$. Suppose that the identity \eqnref{eq:rnhn:obj} holds for some $1\le k<n$. Then, using the recursions \eqnref{eq:recursion:rns} and \eqnref{eq:recursion:hcd} with $m=n-k$, and the induction hypothesis, we have
\begin{align*}
    r_nh_n &= r_{n-k}h_{n-k} + s_{n-k}h_{n-k+1} \\
    &= (c_{n-k}r_{n-k-1}+s_{n-k-1})h_{n-k} + (d_{n-k}r_{n-k-1})h_{n-k+1}\\
    &= r_{n-k-1}(c_{n-k}h_{n-k}+d_{n-k}h_{n-k+1})+s_{n-k-1} h_{n-k}\\
    &=	r_{n-k-1}h_{n-k-1} + s_{n-k-1}h_{n-k}.
\end{align*}
Therefore, the claim \eqnref{eq:rnhn:obj} holds for $k\leq n$, which upon setting $k=n$ and noting the conditions $r_0=1$ and $s_0=0$ give $$r_nh_n=r_0h_0+s_0h_1=h_0.$$ The proof of the identity $\widetilde{r}_n h_n = h_1$ is analogous, and hence it is omitted.
\end{proof}
		
Now if $g$ is nonzero, then the analytic function $h_0(p)=\widehat{g}(p)$ does not have zeros accumulating to $p=0$.  Then by Lemma \ref{lemma:main:coeffrepr}, there exists some $p_1\in(0,p_0)$ satisfying
\beq\label{eq:h:rratio}
	h_1(p)=\frac{\widetilde{r}_n(p)}{r_n(p)}\widehat{g}(p)\quad\mbox{and}\quad h_n(p)=\frac{1}{r_n(p)}\widehat{g}(p),\quad\forall p\in(0,p_1).
\eeq
		
\section{Proof of Theorem \ref{thm:asympt:1stord}} \label{thmsec:asympt:1stord}
In this section, we prove Theorem \ref{thm:asympt:1stord}, using the asymptotic formulas of $h_1(p)$ and $h_n(p)$ as $p\to 0$. These expansions play a key role in the uniqueness proofs in Section \ref{thmsec:Unique:invpb}. The proof of Theorem \ref{thm:asympt:1stord} is fairly lengthy and technical, and it is divided into several steps. More specifically, we proceed in four steps:
\begin{itemize}
    \item[(i)] Derive the asymptotics for the auxiliary functions $E_j(p)$, $F_j(p)$ and $G_j(p)$ as $p\to0$.
    \item[(ii)] Derive the first-order approximations of $h_1(p)$ and $h_n(p)$ via the identities in \eqref{eq:h:rratio}. This requires analyzing $r_n$ and $\tilde r_n$, and motivates the double sequence $X_l^m$ defined in \eqref{eq:def:X}, which enjoys nice properties in Lemmas \ref{lemma:X:property} and \ref{lem:X:positive}, and can express $r_n$ and $\widetilde r_n$, cf. Lemma \ref{lem:rnrtild:asymp}.
    \item[(iii)] Represent the first-order asymptotics of $h_1(p)$ and $h_n(p)$ in terms of solutions to associated boundary value problems, cf. Lemma \ref{lemma:h1hn:approx}.
    \item[(iv)] Prove Theorem \ref{thm:asympt:1stord} by combining Lemma \ref{lem:laptrans:repr} with Lemma \ref{lemma:h1hn:approx}, and properties of the associated boundary value problems (which rely on the construction of the double sequence $X_l^m$ in \eqref{eq:def:X}). 
\end{itemize}
Steps (ii) and (iii) represent the most technical parts of the whole analysis, and the relevant details are given in Section \ref{ssec:asym:rntildern}. These represent the main technical novelty of this section.

\subsection{The first-order approximations of $E_j$, $F_j$ and $G_j$}
We shall approximate the functions $E_j(p)$, $F_j(p)$ and $G_j(p)$ as $p\to0$ by their first-order nonzero terms. We first give an important property of the fundamental solutions $v_j$ and $w_j$.
		
\begin{lemma}\label{lem:v:decreasing}
For all $j=0,1,\dots,n$, we have for all $x^+,x^-\in[x_j,x_{j+1}]\mbox{ satisfying }x^-<x^+$,
\begin{align*}
	\sigma(x^-)v_j'(x^-)\le \sigma(x^+)v_j'(x^+)<0<\sigma(x^-)w_j'(x^-)\le \sigma(x^+)w_j'(x^+).
\end{align*}
\end{lemma}
\begin{proof}
Suppose that $v_j$ is negative at some point in $(x_j,x_{j+1})$.
Then, by the intermediate value theorem, there exists some $a_j\in(x_j,x_{j+1})$ such that $v_j(a_j)=0$.
Solving $-(\sigma v_j')' + q v_j=0$ in $(a_{j},x_{j+1})$ subject to the boundary condition $v_j(a_j)=v_j(x_{j+1})=0$, we have $v_j=0$ in $(a_{j},x_{j+1})$.
Solving $-(\sigma v_j')' + q v_j=0$ in $(x_j,a_{j})$ subject to $v_j(a_j)=v_j'(a_j)=0$, we have $v_j=0$ in $(x_j,a_{j})$, which contradicts the initial condition $v_j(x_j)=1$.
Therefore, we have $v_j\ge0$ for all $j$, which gives
$$(\sigma (x)v_j'(x))' = q(x) v_j(x)\ge0,\quad{\forall }x\in(x_j,x_{j+1}).$$
Thus, $\sigma(x^-)v_j'(x^-)\le \sigma(x^+)v_j'(x^+)$ for all $x^+,x^-\in[x_j,x_{j+1}]$ satisfying $x^-<x^+$.
Next we show $v_j'(x_{j+1})<0$. We have $v_j'(x_{j+1})\le0$ because $v_j\ge0$ and $v_j(x_{j+1})=0$.
If $v_j'(x_{j+1})=0$, solving $-(\sigma v_j')' + q v_j=0$ in $(x_j,x_{j+1})$ subject to $v_j(x_{j+1})=v_j'(x_{j+1})=0$, yields $v_j=0$ in $(x_j,x_{j+1})$, which again contradicts the initial condition $v_j(x_j)=1$.
Thus we arrive at the desired assertions on $v_j$.
The assertions on $w_j$ follow from the above argument for the function $w_j(x_j+x_{j+1}-x)$.
\end{proof}
		
The following lemma is a direct consequence of \eqnref{eq:E:Parseval} and Lemma \ref{lem:v:decreasing}.	\begin{lemma}\label{lem:EFG:signs}
For $j=0,1,\dots,n$, we have
$0<E_{j}^*,G_{j}^*<\infty$ and $  -\infty<F_{j}^*<0.$
\end{lemma}
		
We have the following asymptotic relations.
\begin{lemma}\label{lem:asymp:EFG:0}
For $j=0,1,\dots,n$, as $p\to0$, we have
\begin{align*}
	E_j(p)&= p^{\alpha_j}E_{j}^* + O\left(p^{2\alpha_j}\right),\quad
	F_j(p)= p^{\alpha_j}F_{j}^* + O\left(p^{2\alpha_j}\right),\quad
 G_j(p)= p^{\alpha_j}G_{j}^* + O\left(p^{2\alpha_j}\right).
\end{align*}
\end{lemma}
\begin{proof}
Using the definitions of $E_j(p)$ and $E_j^*$, we have the following geometric series
\begin{align*}
	E_j(p)&
	=p^{\alpha_j}E_j^* - p^{2\alpha_j}\sum_{k=1}^\infty \frac{\big(\sigma(x_j)\varphi_{j,k}'(x_j)\big)^2}{\lambda_{j,k}^3}\sum_{l=0}^\infty (-1)^{l}\left(\frac{p^{\alpha_j}}{\lambda_{j,k}}\right)^{l},\quad\mbox{for }0<p<(\lambda_{j,1})^{1/\alpha_j}.
\end{align*}
Therefore, we have
\begin{align*}
	|E_j(p)-p^{\alpha_j}E_j^*|p^{-2\alpha_j}&\le\sum_{k=1}^\infty \frac{\big(\sigma(x_j)\varphi_{j,k}'(x_j)\big)^2}{\lambda_{j,k}^2\lambda_{j,1}}\sum_{l=0}^\infty \left(\frac{p^{\alpha_j}}{\lambda_{j,1}}\right)^{l}\\			&=\left(\frac{E_j^*}{\lambda_{j,1}}\right)\left(\frac{\lambda_{j,1}}{\lambda_{j,1}-p^{\alpha_j}}\right)\le\frac{2E_j^*}{\lambda_{j,1}},\quad\mbox{for }0<p<(\lambda_{j,1}/2)^{1/\alpha_j}.
\end{align*}
The proofs of the asymptotic relations for $F_j$ and $G_j$ are analogous.
\end{proof}
		
\subsection{The first-order approximation of $h_1$ and $h_n$}\label{ssec:asym:rntildern}
Now we derive the first-order asymptotic behavior of $h_1(p)$ and $h_n(p)$ as $p\to0$. This represents the most technical part of the proof. The challenge lies in the strong coupling among different factors $h_j(p)$. Near $p=0$, $h_1(p)$ and $h_n(p)$ have the expression \eqnref{eq:h:rratio}. Using the definition \eqnref{def:cmdm} and Lemma \ref{lem:asymp:EFG:0}, we obtain that, for $m=1,\dots,n$ and as $p\to 0$,
\begin{align*}
	c_m(p)&=\frac{E_m^*p^{\alpha_m}+G_{m-1}^*p^{\alpha_{m-1}}-\sigma(x_m)v_m'(x_m)+\sigma(x_m)w_{m-1}'(x_{m}) + O(p^{\min\{2\alpha_{m-1},2\alpha_{m}\}})}{F_{m-1}^*p^{\alpha_{m-1}}-\sigma(x_m)v_{m-1}'(x_m)+O(p^{2\alpha_{m-1}})},\\ d_m(p)&=-\frac{F_m^*p^{\alpha_{m}}+\sigma(x_m)w_m'(x_{m})+O(p^{2\alpha_{m}})}{F_{m-1}^*p^{\alpha_{m-1}}-\sigma(x_m)v_{m-1}'(x_m)+O(p^{2\alpha_{m-1}})}.
\end{align*}
For $m=1,\dots,n$, using the identity
\begin{align*}
	&\left(F_{m-1}^*p^{\alpha_{m-1}}-\sigma(x_m)v_{m-1}'(x_m)+O(p^{2\alpha_{m-1}})\right)^{-1}\\
	=&-\frac{1}{\sigma(x_m)v_{m-1}'(x_m)}\left( 1-\frac{F_{m-1}^*}{\sigma(x_m)v_{m-1}'(x_m)}p^{\alpha_{m-1}}+O\left(p^{2\alpha_{m-1}}\right)\right)^{-1}\\
	=&-\frac{1}{\sigma(x_m)v_{m-1}'(x_m)}\left( 1+\frac{F_{m-1}^*}{\sigma(x_m)v_{m-1}'(x_m)}p^{\alpha_{m-1}}+O\left(p^{2\alpha_{m-1}}\right)\right),\quad\mbox{as }p\to0,
\end{align*}
we arrive at the asymptotic formulas
\beq\label{eq:cd:asymptotics}
		\begin{aligned}
			c_m(p)&=c_m^*+ c_{m,0}^*p^{\alpha_m} + c_{m,-}^{*} p^{\alpha_{m-1}} +  O\left(p^{\min\{2\alpha_{m},2\alpha_{m-1}\}}\right),\quad\mbox{as }p\to0,\\
			d_m(p)&= d_m^* + d_{m,0}^* p^{\alpha_m} + d_{m,-}^{*}p^{\alpha_{m-1}}+ O\left(p^{\min\{2\alpha_{m},2\alpha_{m-1}\}}\right),\quad\mbox{as }p\to0,
		\end{aligned}
		\eeq
		where the constants $c_m^*$, $c_{m,0}^*$, $c_{m,-}^*$, $d_m^*$, $d_{m,0}^*$ and $d_{m,-}^*$ are defined by
		\beq\label{eq:cdstar:def}
		\begin{aligned}
			c_m^*:=\frac{v_m'(x_m) - w_{m-1}'(x_{m})}{v_{m-1}'(x_m)},\quad c_{m,0}^*&:= -\frac{E_m^*}{\sigma(x_m)v_{m-1}'(x_m)},\quad c_{m,-}^{*} :=- \frac{G_{m-1}^* - c_m^*F_{m-1}^*}{\sigma(x_m)v_{m-1}'(x_m)},\\
			d_m^*:=\frac{w_m'(x_{m})}{v_{m-1}'(x_m)},\quad
			d_{m,0}^*&:= \frac{F_m^*}{\sigma(x_m)v_{m-1}'(x_m)},\quad d_{m,-}^* := \frac{d_m^* F_{m-1}^*}{\sigma(x_m)v_{m-1}'(x_m)}.
		\end{aligned}
		\eeq

The next lemma gives a crucial inequality for the constants $c_m^*$ and $d_m^*$. 
\begin{lemma}\label{lem:cmdm:bounds}
For all $m=1,\dots,n$, we have
\beq\label{eq:cmdm:bounds}
     c_m^*\ge 1-d_m^*\quad\mbox{and}\quad d_m^*<0.
\eeq
\end{lemma}
\begin{proof}
From Lemma \ref{lem:v:decreasing}, we immediately get $d_m^*<0$. For $m=0,1,\dots,n$, let
\begin{equation*} 
f_m(x):=v_m(x)+w_m(x),\quad  x\in(x_m,x_{m+1}).
\end{equation*}
Then $f_m(x_m)=f_m(x_{m+1})=1$.
We will prove $f_m'(x_m)\le 0\le f_{m-1}'(x_m)$ for all $m=1,\dots,n$, which directly implies $c_m^*\ge 1-d_m^*$ for all $m=1,\dots,n$. Fix $m=0,1,\dots,n$. Then there holds that $f_m>0$ in $(x_m,x_{m+1})$ and
\begin{align*}
    \left\{\begin{aligned}
(\sigma(x)f_m'(x))'+q(x)f_m(x)&=0,\quad \mbox{in }(x_m,x_{m+1}),\\
f_m(x_m)=f_m(x_{m+1})&=1.
\end{aligned}\right.
\end{align*}
If $f_m'(x_m)>0$, then for all $x\in(x_m,x_{m+1})$, we have
\begin{align*}
f_m'(x)&=\frac{\sigma(x)f_m'(x)}{\sigma(x)}=\frac{\sigma(x_m)f_m'(x_m)+\int_{x_m}^x (\sigma(s)f_m(s))'\,{\rm d}s}{\sigma(x)}\\
&=\frac{\sigma(x_m)f_m'(x_m)+\int_{x_m}^x q(s)f_m(s)\,{\rm d}s}{\sigma(x)}>0,
\end{align*}
which gives a contradiction
$$0=f_m(x_{m+1})-f_m(x_m)=\int_{x_m}^{x_{m+1}}f_m'(s)\,{\rm d}s>0.$$
Thus we conclude $f_m'(x_m)\le 0$.
Likewise, if $f_m'(x_{m+1})<0$, then for all $x\in(x_m,x_{m+1})$, there holds
$$f_m'(x)=\frac{\sigma(x)f_m'(x)}{\sigma(x)}=\frac{\sigma(x_{m+1})f_m'(x_{m+1})-\int_x^{x_{m+1}} q(s)f_m(s)\,{\rm d}s}{\sigma(x)}<0,$$
which again gives a contradiction
$$0=f_m(x_{m+1})-f_m(x_m)=\int_{x_m}^{x_{m+1}}f_m'(s)\,{\rm d}s<0.$$
Thus we have $f_{m}'(x_{m+1})\ge 0$. Therefore, $f_m'(x_m)\le 0\le f_{m-1}'(x_m)$ for all $m=1,\dots,n$. This and the definitions of $c_m^*$ and $d_m^*$ give the desired assertion $c_m^*\ge 1-d_m^*$ for all $m=1,\dots,n$.
\end{proof}
		
Next we find the asymptotics of $r_n(p)$ and $\widetilde{r}_n(p)$ as $p\to0$.
It follows from \eqnref{eq:recursion:rns} and \eqnref{eq:cd:asymptotics} that
$$\lim_{p\to0} r_n(p)=X_1^n\quad\mbox{and}\quad \lim_{p\to0} \widetilde{r}_n(p)=X_2^n,$$
where the double sequence $\{X_l^m\}_{l\ge1,\,m\le n}$ is defined by
\beq\label{eq:def:X}
\begin{aligned}
X_l^{l} &:=c_l^*,\quad X_{l}^{l-1}:=1,\quad X_l^{l-i}:=0,\quad\mbox{for all }l=1,2,\dots,n\mbox{ and }i\ge2,\\
X_l^m &:= c_m^* X_l^{m-1} + d_{m-1}^*X_l^{m-2},\quad\mbox{for all }l\mbox{ and }m\mbox{ satisfying }1\le l<m\le n.
\end{aligned}
\eeq
The definition of the double sequence $X_l^m$ follows roughly that of $r_n$ and $\widetilde r_n$, except with the cases $l>2$, which will be needed for representing the asymptotics of $r_n$ and $\widetilde r_n$ in Lemma \ref{lem:rnrtild:asymp} below. See  Lemmas \ref{lemma:X:property} and \ref{lem:X:positive} below for properties of the double sequence $X_l^m$. 

First,  we derive the asymptotic formulas for $r_n(p)$ and $\widetilde{r}_n(p)$ in terms of $X_l^m$. 
		\begin{lemma}\label{lem:rnrtild:asymp}
			For all $m=1,\dots,n$, as $p\to0$, we have
			\begin{align}
			\notag	r_m(p)=X_1^m &+\sum_{i=1}^{m-1}\left[ c_{i+1,-}^* X_1^i X_{i+2}^m + c_{i,0}^* X_1^{i-1} X_{i+1}^m + d_{i+1,-}^* X_1^i X_{i+3}^m + d_{i,0}^* X_1^{i-1} X_{i+2}^m \right] p^{\alpha_{i}}\\
			\label{eq:rnrtild:asymp}	&+ \left(c_{1,-}^* X_2^m+d_{1,-}^* X_3^m\right)p^{\alpha_0} + c_{m,0}^* X_1^{m-1} p^{\alpha_m}
			+O\left(p^{\min\{2\alpha_i\,:\,0\le i\le m\}}\right),\\
				\notag	\widetilde{r}_m(p)=X_2^m &+\sum_{i=1}^{m-1}\left[ c_{i+1,-}^* X_2^i X_{i+2}^m + c_{i,0}^* X_2^{i-1} X_{i+1}^m + d_{i+1,-}^* X_2^i X_{i+3}^m + d_{i,0}^* X_2^{i-1} X_{i+2}^m \right] p^{\alpha_{i}} \\
			\label{eq:rnrtild:asymp:tild}	&+ c_{m,0}^* X_2^{m-1} p^{\alpha_m}+O\left(p^{\min\{2\alpha_i\,:\,1\le i\le m\}}\right).
			\end{align}
		\end{lemma}
\begin{proof}
The proof is based on mathematical induction on $m$. From the recursions \eqnref{eq:recursion:rns} and \eqnref{eq:def:X} and the approximation \eqnref{eq:cd:asymptotics}, we have, as $p\to0$,
\begin{align*}
\widetilde{r}_1(p)&=1=X_2^1,\\
r_1(p)&=c_1(p)=X_1^1 +c_{1,0}^*  p^{\alpha_1} +c_{1,-}^* p^{\alpha_0} +O\left(p^{\min\{2\alpha_0,2\alpha_1\}}\right).
\end{align*}
Thus, the identities \eqref{eq:rnrtild:asymp} and \eqref{eq:rnrtild:asymp:tild} hold for $m=1$.
Likewise, for $m=2$, we have    \begin{align*}&
\widetilde{r}_2(p)=c_2(p)=X_2^2 + c_{2,0}^* p^{\alpha_2} +c_{2,-}^* p^{\alpha_1} +O\left(p^{\min\{2\alpha_1,2\alpha_2\}}\right),
\end{align*}
and similarly, we have
\begin{align*}
r_2(p)=&c_1(p)c_2(p)+d_1(p)				=c_1^*c_2^*+d_1^*+ d_{1,0}^*p^{\alpha_1} + d_{1,-}^{*}p^{\alpha_{0}}\\
 &+ c_2^*\left(c_{1,0}^*p^{\alpha_1} + c_{1,-}^{*} p^{\alpha_{0}} \right) + c_1^*\left( c_{2,0}^*p^{\alpha_2} + c_{2,-}^{*} p^{\alpha_{1}} \right) +O\left(p^{\min\{2\alpha_i\,:\,i\le2\}}\right)\\
=&X_1^2  + \left( c_{2,-}^* X_1^1 +  c_{1,0}^* X_2^2 +d_{1,0}^*\right)p^{\alpha_1} + \left(c_{1,-}^* X_2^2 +d_{1,-}^*\right)p^{\alpha_0} \\
 &+ c_{2,0}^* X_1^1 p^{\alpha_2}+O\left(p^{\min\{2\alpha_i\,:\,i\le2\}}\right).
\end{align*}
Thus, the identities \eqnref{eq:rnrtild:asymp} and \eqnref{eq:rnrtild:asymp:tild} hold for $m=2$.
Next we prove \eqnref{eq:rnrtild:asymp} using mathematical induction. Suppose that \eqnref{eq:rnrtild:asymp} holds for $m=1,\dots,k$, with $2\le k<n$.
By the  recursion \eqnref{eq:recursion:rns} and the induction hypothesis (i.e., \eqnref{eq:rnrtild:asymp} with $m=k$ and $m=k-1$), we obtain
\begin{align}
		&r_{k+1}(p)=c_{k+1}(p) r_k(p) + d_{k}(p) r_{k-1}(p)\nonumber\\
		=&c_{k+1}(p) X_1^k + d_{k}(p) X_1^{k-1} + R_k(p)+ O\left(p^{\min\{2\alpha_i\,:\,0\le i\le k\}}\right),\quad\mbox{as }p\to0,\label{eq:r:intermsof:X}
    \end{align}
where the term $R_k$ is given by $R_k=R_k^{(1)}+R_k^{(2)}$, with $R_k^{(1)}$ and $R_k^{(2)}$ defined respectively by
\begin{align*}
    R_k^{(1)}(p)&:=\sum_{i=1}^{k-2}\left[c_{i+1,-}^* X_1^i (c_{k+1}^* X_{i+2}^{k} +d_{k}^* X_{i+2}^{k-1} ) + c_{i,0}^* X_1^{i-1} (c_{k+1}^* X_{i+1}^{k} +d_{k}^* X_{i+1}^{k-1} ) \right]p^{\alpha_i}\\
	&\quad +\sum_{i=1}^{k-2}\left[ d_{i+1,-}^* X_1^i (c_{k+1}^* X_{i+3}^{k} +d_{k}^* X_{i+3}^{k-1} ) + d_{i,0}^* X_1^{i-1} (c_{k+1}^* X_{i+2}^{k} +d_{k}^* X_{i+2}^{k-1} )\right]p^{\alpha_i},\\
			R_k^{(2)}(p)&:= c_{k+1}^* \left(c_{k,-}^* X_1^{k-1} + c_{k-1,0}^* X_1^{k-2}  X_{k}^{k}  + d_{{k-1},0}^* X_1^{k-2} \right) p^{\alpha_{k-1}}\\
		&\quad + \left[c_{1,-}^* (c_{k+1}^*X_2^k + d_{k}^* X_2^{k-1}) + d_{1,-}^*(c_{k+1}^*X_3^k + d_{k}^* X_3^{k-1})\right]p^{\alpha_0}\\
	& \quad + c_{k+1}^* \left(c_{k,0}^*X_1^{k-1}p^{\alpha_k}\right)+ d_{k}^* \left(c_{k-1,0}^*X_1^{k-2}p^{\alpha_{k-1}}\right).
\end{align*}
It remains to simplify the expressions of $R_k^{(1)}$ and $R_k^{(2)}$. Using the defining recursion \eqnref{eq:def:X} of $\{X_l^m\}_{l\ge1,\,m\le n}$, we get
		\begin{align}
			\label{eq:Rk1:simp}
			R_k^{(1)}(p)&=\sum_{i=1}^{k-2}\left[c_{i+1,-}^* X_1^i  X_{i+2}^{k+1} + c_{i,0}^* X_1^{i-1}  X_{i+1}^{k+1} + d_{i+1,-}^* X_1^i  X_{i+3}^{k+1} + d_{i,0}^* X_1^{i-1} X_{i+2}^{k+1} \right]p^{\alpha_i},\\
			\notag R_k^{(2)}(p)&= \left[c_{k+1}^* \left(c_{k,-}^* X_1^{k-1} + d_{{k-1},0}^* X_1^{k-2} \right)+ c_{k-1,0}^* X_1^{k-2} X_k^{k+1}\right]p^{\alpha_{k-1}}\\
			&\quad + \left(c_{1,-}^* X_2^{k+1} + d_{1,-}^* X_3^{k+1}\right)p^{\alpha_0} +  c_{k+1}^* c_{k,0}^*X_1^{k-1}p^{\alpha_k}.
		\end{align}
It follows from the identity \eqnref{eq:cd:asymptotics} and the recursion \eqnref{eq:def:X} that we also have, as $p\to0$,
		\begin{align*}
			c_{k+1}(p) X_1^k + d_{k}(p) X_1^{k-1}&=X_1^{k+1}+\left(c_{k+1,0}^*  p^{\alpha_{k+1}} +c_{k+1,-}^* p^{\alpha_k} \right) X_1^k \\
			&\quad + \left(d_{k,0}^*  p^{\alpha_{k}} + d_{k,-}^* p^{\alpha_{k-1}} \right) X_1^{k-1} + O\left(p^{\min\{2\alpha_{k-1},2\alpha_{k},2\alpha_{k+1}\}}\right).
		\end{align*}
Therefore, we have, as $p\to0$,
		\beq\label{eq:aux:cpdpR}
		\begin{aligned}
			&c_{k+1}(p) X_1^k + d_{k}(p) X_1^{k-1} + R_k^{(2)}(p)\\
			=&X_1^{k+1}+\left(c_{k+1,-}^* X_1^k + c_{k+1}^*c_{k,0}^*X_1^{k-1}  + d_{k,0}^*X_1^{k-1}\right)p^{\alpha_k}\\
			&+\left[c_{k+1}^* \left(c_{k,-}^* X_1^{k-1} + d_{{k-1},0}^* X_1^{k-2} \right)+ c_{k-1,0}^* X_1^{k-2} X_k^{k+1} + d_{k,-}^* X_1^{k-1}\right]p^{\alpha_{k-1}}\\
			&+ \left(c_{1,-}^* X_2^{k+1} + d_{1,-}^* X_3^{k+1}\right)p^{\alpha_0}+c_{k+1,0}^*  X_1^k p^{\alpha_{k+1}} + O\left(p^{\min\{2\alpha_{k-1},2\alpha_{k},2\alpha_{k+1}\}}\right).
		\end{aligned}
	\eeq
Further, from the definition \eqnref{eq:def:X} again, we have
		\beq\label{eq:ieqkorkm1:summand}
		\begin{aligned}
		&c_{i+1,-}^* X_1^i  X_{i+2}^{k+1} + c_{i,0}^* X_1^{i-1}  X_{i+1}^{k+1} + d_{i+1,-}^* X_1^i  X_{i+3}^{k+1} + d_{i,0}^* X_1^{i-1} X_{i+2}^{k+1}\\
		=&\begin{cases}
			\displaystyle c_{k+1}^* \left(c_{k,-}^* X_1^{k-1} + d_{{k-1},0}^* X_1^{k-2} \right)+ c_{k-1,0}^* X_1^{k-2} X_k^{k+1} + d_{k,-}^* X_1^{k-1},&\mbox{if }i=k-1,\\
			\displaystyle c_{k+1,-}^* X_1^k + c_{k+1}^*c_{k,0}^*X_1^{k-1}  + d_{k,0}^*X_1^{k-1},&\mbox{if }i=k.
		\end{cases}
		\end{aligned}
	\eeq
Substituting \eqnref{eq:Rk1:simp}, \eqnref{eq:aux:cpdpR} and \eqnref{eq:ieqkorkm1:summand} into \eqnref{eq:r:intermsof:X} yields the identity \eqnref{eq:rnrtild:asymp} with $m=k+1$, i.e., the induction step, which completes the proof of the identity \eqnref{eq:rnrtild:asymp} for all $m=1,\dots,n$.
The assertion \eqnref{eq:rnrtild:asymp:tild} can be proved similarly, and hence the details are omitted.
\end{proof}
		
The double sequence $X_m^n$ satisfies the following property.
\begin{lemma}\label{lemma:X:property}
For $m=1,\dots,n$, we have
\begin{equation*}
   c_m^* X_{m+1}^n + d_m^* X_{m+2}^n = X_{m}^n
	\quad\mbox{and}\quad
	X_2^m X_1^n-X_1^m X_2^n=(-1)^{m+1}X_{m+2}^{n}\prod_{j=1}^{m} \frac{w_j'(x_{j})}{v_{j-1}'(x_{j})}.
\end{equation*}
\end{lemma}
\begin{proof}
By mathematical induction on $k$, we will prove the assertion
\beq\label{eq:cxpdx:subscript}
	c_m^* X_{m+1}^{m+k} + d_m^* X_{m+2}^{m+k} = X_{m}^{m+k},\quad\mbox{for all }1\le m\le n\mbox{ and }0\le k\le n-m.
\eeq
From the definition \eqnref{eq:def:X}, we have, for all $m=1,\dots,n$,
\begin{align*}
		&X_m^m=c_m^*,\quad X_{m+1}^m=1,\quad X_{m+2}^m=0,\\	&X_{m}^{m+1}=c_{m+1}^*X_{m}^{m}+d_{m}^* X_{m}^{m-1}=c_{m+1}^* c_{m}^*+d_{m}^*,
\end{align*}
which gives the claim \eqnref{eq:cxpdx:subscript} with $k=0$ and $k=1$.
This also proves \eqnref{eq:cxpdx:subscript} with $n-m\le1$, so we may assume $n-m\ge 2$ below.
Now suppose that \eqnref{eq:cxpdx:subscript} holds for all $m=1,\dots,n$ and  $0\le k<n-m$.
Then, using \eqnref{eq:def:X} and the induction hypothesis, we obtain
\begin{align*}
	&X_{m}^{m+k+1}=c_{m+k+1}^*X_{m}^{m+k}+d_{m+k}^* X_{m}^{m+k-1}\\
	=&c_{m+k+1}^*\left(c_m^* X_{m+1}^{m+k} + d_m^* X_{m+2}^{m+k}\right)+d_{m+k}^* \left(c_m^* X_{m+1}^{m+k-1} + d_m^* X_{m+2}^{m+k-1}\right).
\end{align*}
Changing the order of summation and using the definition \eqnref{eq:def:X} again, we arrive at
\begin{align*}
	X_{m}^{m+k+1}&=c_m^*\left(c_{m+k+1}^* X_{m+1}^{m+k} + d_{m+k}^*X_{m+1}^{m+k-1}\right) + d_m^*\left(c_{m+k+1}^* X_{m+2}^{m+k}+ d_{m+k}^* X_{m+2}^{m+k-1}\right)\\
	&=c_m^* X_{m+1}^{m+k+1} + d_m^* X_{m+2}^{m+k+1},
\end{align*}
which proves \eqnref{eq:cxpdx:subscript}.
This gives the first identity by taking $k=n-m$.
To prove the second identity,
using the first identity repeatedly, we derive
\begin{align*}
X_2^m X_1^n-X_1^m X_2^n &= X_2^m\left(c_1^*X_2^n+d_1^*X_3^n\right) - \left(c_1^*X_2^m+d_1^*X_3^m\right)X_2^n\\
&= d_1^*\left(X_2^m X_3^n-X_3^m X_2^n\right)= d_1^*d_2^*\left(X_4^m X_3^n-X_3^m X_4^n\right) \\
&= \cdots = \bigg(\prod_{j=1}^{m} d_j^*\bigg)(-1)^{m+1}\left(X_{m+1}^m X_{m+2}^n-X_{m+2}^m X_{m+1}^n\right).
\end{align*}
Since $X_{m+1}^m=1$ and $X_{m+2}^m=0$, the identity $\prod_{j=1}^{m} d_j^*=\prod_{j=1}^{m} \frac{w_j'(x_{j})}{v_{j-1}'(x_{j})}$
gives the second assertion.
\end{proof}
		
The next result gives the positivity of $X_m^n$. The proof is based on \eqnref{eq:cmdm:bounds}.
\begin{lemma}\label{lem:X:positive}
We have $X_m^n>0$ for $m=1,2,\dots,n$.
\end{lemma}
\begin{proof}
Fix $m=1,2,\dots,n$ and consider the sequences $\{y_i\}_{i=m-1}^n$ and $\{z_i\}_{i=m-1}^n$ defined by
\begin{align*}
	&y_{m-1}=1,\quad z_{m-1}=0,\\
&y_i=c_i^* y_{i-1} + z_{i-1},\quad z_i=d_i^* y_{i-1},\quad\mbox{for }i=m,m+1,\dots,n.
\end{align*}
Then we have $y_n=X_m^n$, cf. the definition \eqref{eq:def:X}.
We will prove
\beq\label{eq:Xpos:obj}
  d_{i}^*y_{i-1}+y_i\ge 1,\quad\mbox{for }i=m,m+1,\dots,n,
\eeq
by mathematical induction (on $i$). From $y_{m}=c_1^*$ and \eqnref{eq:cmdm:bounds} with $i=m$, we have
\beq\label{eq:dypy:1}
	d_m^*y_{m-1}+y_m=d_m^*+c_m^*\ge d_m^*+(1-d_m^*)=1.
\eeq
That is, the claim \eqref{eq:Xpos:obj} holds for $i=m$. Next, suppose $d_{i}^*y_{i-1}+y_i\ge1$ for all $i$ satisfying $m\le i\le k< n$. Then, by Lemma \ref{lem:cmdm:bounds}, we have 
			$$y_k>(-d_k^*)y_{k-1}>\cdots>\bigg(\prod_{j=m}^k (-d_j^* )\bigg)y_{m-1}>0.$$
			Using this inequality, \eqnref{eq:cmdm:bounds} and the induction hypothesis $d_{k}^*y_{k-1}+y_k\ge 1$, we arrive at
\begin{align}
d_{k+1}^*y_{k}+y_{k+1} &= d_{k+1}^*y_{k}+\left(c_{k+1}^* y_k + d_k^*y_{k-1}\right)\nonumber	\\
&\ge d_{k+1}^*y_{k}+\left((1-d_{k+1}^*) y_k + d_k^*y_{k-1}\right)=d_{k}^*y_{k-1}+y_k\ge 1.
\label{eq:dypy:k}			\end{align}
Hence, we have proved the claim \eqnref{eq:Xpos:obj}.
Finally, using the assertion \eqnref{eq:Xpos:obj}, we obtain
\begin{equation*}
X_m^n=y_n>(-d_n^*)y_{n-1}>(-d_n^*)(-d_{n-1}^*)y_{n-2}>\cdots>\bigg(\prod_{j=m}^n (-d_j^* )\bigg)y_{m-1}>0,
\end{equation*}
which completes the proof of the lemma.
\end{proof}

To succinctly represent the factors $h_1(p)$ and $h_n(p)$, now we define the following three functions $u,\overline{u},\widetilde{u}:[0,L]\to \mathbb{R}$, for all $i=0,1,\dots,n$ and $x\in(x_i,x_{i+1})$, by
    \begin{align}
        \label{u:def:exp}
		u(x)&:=\frac{X_{i+1}^n}{X_1^n}v_i(x) + \frac{X_{i+2}^n}{X_1^n}w_i(x),\\
        \label{ubar:def:exp}
            \overline{u}(x)&:=\frac{X_1^{i-1}}{X_1^n}\frac{v_n'(L)}{v_{i-1}'(x_i)}\frac{\sigma(L)}{\sigma(x_i)}v_i(x)+\frac{X_1^{i}}{X_1^n}\frac{v_n'(L)}{v_{i}'(x_{i+1})}\frac{\sigma(L)}{\sigma(x_{i+1})}w_i(x),\\
        \label{v:def:exp}
		\widetilde{u}(x)&:=\frac{X_{i+1}^n}{X_1^n}M_{i-1}v_i(x) + \frac{X_{i+2}^n}{X_1^n}M_iw_i(x),
    \end{align}
where the sequence $\{M_i\}_{i=0}^n$ {\color{blue}is} defined by
\begin{equation}\label{eq:M:def}
    M_{-1}:=1\quad \mbox{and}\quad M_{i}:=\frac{\sigma(0)}{\sigma(x_{i+1})}\prod_{j=0}^{i}\frac{w_j'(x_{j})}{-v_j'(x_{j+1})},\quad i=0,1,\dots,n.
\end{equation}
    
The properties of these three functions will be discussed in Lemmas \ref{lemma:u:ivp:C0dep} and \ref{lem:X-repre-baru}. Using the functions $u$, $\bar u$ and $\widetilde u$, we finally obtain a concise asymptotic formulas of $h_1(p)$ and $h_n(p)$. This is the main result of this part.
\begin{lemma}\label{lemma:h1hn:approx}
The functions $h_1(p)$ and $h_n(p)$ have the following asymptotics as $p\to0$:
\begin{align}
        \notag
    &h_1(p)\\
    =&\widehat{g}(p)\bigg(\!\!\frac{X_2^n}{X_1^n}-\frac{1}{\sigma(0)w_0'(0)}\bigg [p^{\alpha_0}(u,\widetilde u\!-\!v_0)_{L_\rho^2( 0,x_1)}+\sum_{i=1}^n p^{\alpha_i}(u,\widetilde u)_{L^2_\rho(x_i,x_{i+1})}\bigg] + O\left(p^{2\min_i\{\alpha_i\}}\right)\!\!\bigg), \label{eq:h1:asympt}\\      
  &\notag h_n(p)\\
   =&\widehat{g}(p)\bigg(\!\!\frac{1}{X_1^n}\! +\! \frac{1}{\sigma(L)v_n'(L)}\bigg[\sum_{i=0}^{n\!-\!1}p^{\alpha_i}(u,\overline{u})_{L^2_\rho (x_i,x_{i+1})}\! +\!p^{\alpha_n}(u,\overline{u}\!-\!w_n)_{L^2_\rho(x_n,L)}\bigg]\!+\! O\left(p^{2\min_i\{\alpha_i\}}\right)\!\!\bigg).
            \label{eq:hn:asympt}
\end{align}
\end{lemma}
\begin{proof}
From \eqnref{eq:cdstar:def} and Lemma \ref{lemma:X:property}, we have, for all $0\le i\le n-1$ and $1\le j\le n$,
\begin{align}
		c_{i+1,-}^* X_{i+2}^n + d_{i+1,-}^* X_{i+3}^n
        &=-\frac{G_i^*-c_{i+1}^*F_i^*}{\sigma(x_{i+1})v'_i(x_{i+1})}X_{i+2}^n+\frac{d_{i+1}^*F_i^*}{\sigma(x_{i+1})v'_i(x_{i+1})}X_{i+3}^n\nonumber\\
        &=-\frac{G_i^*}{\sigma(x_{i+1})v'_i(x_{i+1})}X_{i+2}^n+\frac{F_i^*}{\sigma(x_{i+1})v'_i(x_{i+1})}\left(c_{i+1}^*X_{i+2}^n+d_{i+1}^*X_{i+3}^n\right)\nonumber\\
        &=-\frac{G_i^*}{\sigma(x_{i+1})v_i'(x_{i+1})}X_{i+2}^n + \frac{F_i^*}{\sigma(x_{i+1})v_i'(x_{i+1})}X_{i+1}^n,\label{eq:cxpdx:EFG-1}\\
		c_{j,0}^* X_{j+1}^n + d_{j,0}^* X_{j+2}^n&= -\frac{E_j^*}{\sigma(x_{j})v_{j-1}'(x_j)}X_{j+1}^n + \frac{F_j^*}{\sigma(x_{j})v_{j-1}'(x_j)}X_{j+2}^n.\label{eq:cxpdx:EFG-2}
\end{align}
Using Lemma \ref{lem:rnrtild:asymp} with $m=n$ and the identities \eqnref{eq:cxpdx:EFG-1} and \eqref{eq:cxpdx:EFG-2}, we obtain as $p\to0$:
		\begin{align*}
			&\frac{1}{r_n(p)}-\frac{1}{X_1^n}+ O(p^{2\min_i\{\alpha_i\}})\\
			=&-\frac{1}{(X_1^n)^2} \bigg[ \sum_{i=0}^{n-1}\left( c_{i+1,-}^* X_1^i X_{i+2}^n + c_{i,0}^* X_1^{i-1} X_{i+1}^n + d_{i+1,-}^* X_1^i X_{i+3}^n + d_{i,0}^* X_1^{i-1} X_{i+2}^n \right) p^{\alpha_{i}}\\
      &+ c_{n,0}^* X_{1}^{n-1} p^{\alpha_n}\bigg]\\
			=&-\frac{1}{(X_1^n)^2}  \bigg[\sum_{i=0}^{n-1}  \frac{p^{\alpha_{i}}X_1^i}{\sigma(x_{i+1})v_i'(x_{i+1})}\left(-G_i^* X_{i+2}^n + F_i^* X_{i+1}^n\right)\\
    &+ \sum_{i=0}^{n}\frac{p^{\alpha_{i}}X_1^{i-1}}{\sigma(x_i)v_{i-1}'(x_i)}\left(-E_i^* X_{i+1}^n + F_i^* X_{i+2}^n\right)  \bigg].
		\end{align*}
Then, from the representations in \eqnref{eq:E:Parseval}, we obtain as $p\to 0$,
\begin{align*}
    \frac{1}{r_n(p)}=\frac{1}{X_1^n} + \frac{1}{\sigma(L)v_n'(L)}\left[\sum_{i=0}^n p^{\alpha_i}(u,\overline{u})_{L^2_\rho(x_i,x_{i+1})} - p^{\alpha_n}(u,w_n)_{L_\rho^2(x_n,L)}\right]+O\left(p^{2\min_i\{\alpha_i\}}\right).
\end{align*}
Combining this identity with \eqnref{eq:h:rratio} yields the identity \eqnref{eq:hn:asympt}. Next, we prove the assertion \eqnref{eq:h1:asympt}. It follows from Lemma \ref{lem:rnrtild:asymp} with $m=n$ that
		\begin{align*}
			&\frac{\widetilde{r}_n}{r_n}-\frac{X_2^n}{X_1^n} + O\left(p^{2\min_i\{\alpha_i\}}\right)\\
			=&
			-\frac{X_2^n}{(X_1^n)^2}\sum_{i=0}^n\left[ c_{i+1,-}^* X_1^i X_{i+2}^n + c_{i,0}^* X_1^{i-1} X_{i+1}^n + d_{i+1,-}^* X_1^i X_{i+3}^n + d_{i,0}^* X_1^{i-1} X_{i+2}^n \right] p^{\alpha_{i}}\\
			&+	\frac{X_1^n}{(X_1^n)^2}\sum_{i=0}^n\left[ c_{i+1,-}^* X_2^i X_{i+2}^n + c_{i,0}^* X_2^{i-1} X_{i+1}^n + d_{i+1,-}^* X_2^i X_{i+3}^n + d_{i,0}^* X_2^{i-1} X_{i+2}^n \right] p^{\alpha_{i}}.
		\end{align*}
		From \eqnref{eq:h:rratio} and the second identity of Lemma \ref{lemma:X:property}, there holds that as $p\to0$,
\begin{align}
	\frac{h_1(p)}{\widehat{g}(p)}&=\frac{\widetilde{r}_n}{r_n}=\frac{X_2^n}{X_1^n}-p^{\alpha_0}\frac{X_2^n}{(X_1^n)^2}\left(c_{1,-}^* X_2^n + d_{1,-}^* X_3^n\right)\nonumber\\
	&\quad
	-\sum_{i=1}^n p^{\alpha_i} \frac{X_{i+2}^n}{(X_1^n)^2}\left(c_{i+1,-}^* X_{i+2}^n + d_{i+1,-}^* X_{i+3}^n\right)\prod_{j=1}^i\frac{w_j'(x_j)}{-v_{j-1}'(x_j)}\nonumber\\
    &\quad -\sum_{i=1}^n p^{\alpha_i}\frac{X_{i+1}^n}{(X_1^n)^2}\left(c_{i,0}^* X_{i+1}^n + d_{i,0}^* X_{i+2}^n\right)\prod_{j=1}^{i-1}\frac{w_j'(x_j)}{-v_{j-1}'(x_j)} + O\left(p^{2\min_i\{\alpha_i\}}\right).\label{eq:h1:asymp:last}
\end{align}
Combining the identities \eqnref{eq:cxpdx:EFG-1}--\eqref{eq:cxpdx:EFG-2}, \eqnref{eq:E:Parseval} and the definitions \eqnref{u:def:exp} and \eqnref{v:def:exp} yields \eqnref{eq:h1:asympt}. 
\end{proof}

\subsection{Proof of Theorem \ref{thm:asympt:1stord}}
\label{subsect:C:asympt}
Now we prove Theorem \ref{thm:asympt:1stord}, i.e., asymptotics of $\p_x\widehat{U}(p,0)$ and $\p_x\widehat{U}(p,L)$ as $p\to0$. We discuss the two cases separately. 
First we derive the expression for $\p_x\widehat{U}(p,0)$. It follows from Lemma \ref{lem:laptrans:repr} that
		\beq\label{eq:dhatU:at0}
		\begin{aligned}
			\sigma(0)\p_x\widehat{U}(p,0)=\widehat{g}(p)\sigma(0)v_0'(0)-h_{1}(p)\sigma(0)w_0'(0) - \widehat{g}(p)E_0(p) + h_{1}(p)F_0(p).
		\end{aligned}
		\eeq
Using the identity \eqnref{eq:h1:asympt}, we obtain
		\begin{align*}
			\p_x\widehat{U}(p,0)&=\frac{1}{\sigma(0)}\widehat{g}(p)\left(\sigma(0)v_0'(0)-E_0^*p^{\alpha_0} - h_1(p)(\sigma(0)w_0'(0)-F_0^* p^{\alpha_0})+O\left(p^{2\alpha_0}\right)\right)\\
			&=\widehat{g}(p)\left[C_0 + \sum_{i=0}^n C_{i+1} p^{\alpha_i} + O\left(p^{2\min_i\{\alpha_i\}}\right)\right], \quad\mbox{as }p\to0,
		\end{align*}
where the constants $C_m$ are defined by	\beq\label{eq:C0C1C2:def}
    \left\{\begin{aligned}
		C_0&=v_0'(0)+\frac{X_2^n}{X_1^n}w_0'(0),\\
			C_{m+1}&=-(u,\widetilde{u})_{L^2_\rho(x_m,x_{m+1})}/\sigma(0), \quad  m=0,1,\dots,n.
\end{aligned}\right.
\eeq
Meanwhile, the function $u$ defined in \eqref{u:def:exp} has the following property.
\begin{lemma}\label{lemma:u:ivp:C0dep}
The function $u$ defined in \eqref{u:def:exp} belongs to $C^{1,1}[0,L]$ and satisfies 
\beq\label{eq:ivp:u:C0dep}
   \left\{\begin{aligned}
		-(\sigma(x)u'(x))' + q(x) u(x)&=0,\quad\mbox{in }(0,L),\\
      u(0)=1,\quad u'(0)&=C_0.
		\end{aligned}\right.
    \eeq   
\end{lemma}
\begin{proof}
We first recall a gluing property for Lipschitz continuity: for any $a<b<c$ and $f:[a,c]\to\mathbb{R}$, if $f$ and $f'$ are continuous and $f'$ is Lipschitz on both $[a,b]$ and $[b,c]$, then $f\in C^{1,1}[a,c]$.
Indeed, let $L_a$ and $L_b$ be the Lipschitz constants of $f'$ on $[a,b]$ and $[b,c]$, respectively. Then we claim
\beq\label{eq:Lips:gluing}
    |f'(x)-f'(y)|\le \max(L_a,L_b)|x-y|,\quad\forall x,y\in[a,c].
\eeq
Fix any $x,y\in[a,c]$, with $x\le y$. If $x,y\in[a,b]$ or $x,y\in[b,c]$, the claim \eqnref{eq:Lips:gluing} holds trivially. Otherwise, for $x\in[a,b]$ and $y\in[b,c]$, we have
    $$|f'(x)-f'(y)|\le|f'(x)-f'(b)|+|f'(b)-f'(y)|\le L_a(b-x) + L_b(y-b)\le \max(L_a,L_b)(y-x).$$
This proves the claim \eqnref{eq:Lips:gluing}, and the gluing property.  For each $i=0,1,\dots,n$, the functions $v_i$ and $w_i$ belong to $C^{1,1}[x_i,x_{i+1}]$, thus $u\in C^{1,1}[x_i,x_{i+1}]$.
    For each $i=1,\dots,n$, 
    $$\lim_{s\to0^+}u(x_i+s)=\frac{X_{i+1}^n}{X_1^n}=\lim_{s\to0^+}u(x_i-s).$$
    Thus, $u$ is continuous at $x_i$.
    For each $i=1,\dots,n$,  we have
    \begin{align*}
        &\lim_{s\to0^+}\left(\p_x u(x_i + s)-\p_x u(x_i - s)\right)
        =\frac{X_{i+2}^n}{X_1^n}w_i'(x_i) + \frac{X_{i+1}^n}{X_1^n}\left(v_i'(x_i) - w_{i-1}'(x_i)\right)-\frac{X_{i}^n}{X_1^n} v_{i-1}'(x_i).
\end{align*}
From the identity $X_i^n=c_i^*X_{i+1}^n+d_i^*X_{i+2}^n$ from Lemma \ref{lemma:X:property} and the definitions of $c_i^*$ and $d_i^*$, we deduce
\begin{align*}
&\lim_{s\to0^+}\left(\p_x u(x_i + s)-\p_x u(x_i - s)\right)\\
        =&\frac{X_{i+2}^n}{X_1^n}w_i'(x_i) + \frac{X_{i+1}^n}{X_1^n}\left(v_i'(x_i) - w_{i-1}'(x_i)\right)-\frac{(c_i^*X_{i+1}^n+d_i^*X_{i+2}^n)}{X_1^n} v_{i-1}'(x_i)\\
        =&\frac{X_{i+2}^n}{X_1^n}\left(w_i'(x_i)-d_i^* v_{i-1}'(x_i)\right) + \frac{X_{i+1}^n}{X_1^n}\left(v_i'(x_i) - w_{i-1}'(x_i)-c_i^* v_{i-1}'(x_i)\right)=0.
    \end{align*}
Thus, $u'$ is continuous at $x_i$.
By the gluing property, we deduce $u\in C^{1,1}[0,L]$. The relations in \eqnref{eq:ivp:u:C0dep} follow directly from the definitions of $u$ and $C_0$.
\end{proof}

The function $\widetilde{u}$ defined  in \eqnref{v:def:exp} satisfies
\begin{align*}
   \left\{\begin{aligned}
	(\sigma(x)\widetilde{u}'(x))' + q(x) \widetilde{u}(x)&=0,\quad\mbox{in }(x_m,x_{m+1}),\\
	\widetilde{u}(x_m)&=M_{m-1}u(x_m),\\ \widetilde{u}(x_{m+1})&=M_{m}u(x_{m+1}).
	\end{aligned}\right.
\end{align*}

The next lemma shows $u=\widetilde{u}$ in $(0,L)$, which completes the proof of Theorem \ref{thm:asympt:1stord} on $\partial_x\widehat{U}(p,0)$.
\begin{lemma}\label{lemma:Meq1}
For the constants $M_m$ defined in \eqref{eq:M:def}, there holds
$M_m=1$, $m=0,1,\dots,n.$
\end{lemma}
\begin{proof}
Fix $i=0,1,\dots,n$ and define $$g(x):=\sigma(x)\left({v}_i'(x){w}_i(x)-{v}_i(x){w}_i'(x)\right),\quad\mbox{for all }x\in[x_i,{x}_{i+1}].$$
Then $g$ is Lipschitz continuous in $[x_i,{x}_{i+1}]$ and satisfies
\beq\label{eq:g:bdycond}
    g(x_i)=-\sigma(x_i){w}_i'(x_i)\quad\mbox{and}\quad g({x}_{i+1})=\sigma({x}_{i+1}){v}_i'({x}_{i+1}).
\eeq
We also have, for a.e. $x\in[x_i,{x}_{i+1}]$,
\begin{align}
   \notag
    g'(x)&=\sigma(x)\left({v}_i'(x){w}_i'(x)-{v}_i'(x){w}_i'(x)\right) + \left(\sigma(x){v}_i'(x)\right)'{w}_i(x) - \left(\sigma(x){w}_i'(x)\right)'{v}_i(x)\\
   &=q(x){v}_i(x){w}_i(x)-q(x){w}_i(x){v}_i(x)=0.
 \label{eq:g:derivative}        
\end{align}
From \eqnref{eq:g:bdycond} and \eqnref{eq:g:derivative}, we obtain $$-\sigma(x_i){w}_i'(x_i)=g(x_i)=g(x_{i+1})=\sigma(x_{i+1})v_i'(x_{i+1}).$$
Therefore, using the definition \eqnref{eq:M:def}, we arrive at
$$M_m=\prod_{i=0}^{m}\frac{\sigma(x_i)w_i'(x_{i})}{-\sigma(x_{i+1})v_i'(x_{i+1})}=\prod_{i=0}^{m}1=1,\quad m=0,1,\dots,n.$$
This completes the proof of the lemma. 
\end{proof}
		
Now we derive the asymptotic formula of  $\p_x\widehat{U}(p,L)$ as $p\to0$.
It follows from Lemma \ref{lem:laptrans:repr} that
\beq\label{eq:dhatU:atL}
		\begin{aligned}
			\sigma(L)\p_x\widehat{U}(p,L)&=h_n(p)\left(\sigma(L) v_n'(L)-F_n(p)\right).
		\end{aligned}
		\eeq
Applying Lemmas \ref{lem:asymp:EFG:0} and \ref{lemma:h1hn:approx} yields
\begin{align*}
	&\sigma(L)\p_x\widehat{U}(p,L)			=h_n(p)\left(\sigma(L) v_n'(L) - p^{\alpha_n} F_n^* + O\left(p^{2\alpha_n}\right)\right)\\
	=&\widehat{g}(p)\left[\frac{\sigma(L) v_n'(L)}{X_1^n} -p^{\alpha_n}\frac{F_n^*}{X_1^n} + \sum_{i=0}^n p^{\alpha_i}(u,\overline{u})_{L^2_\rho(x_i,x_{i+1})} + O\left(p^{2\min_i\{\alpha_i\}}\right)\right].
\end{align*}
Therefore, we arrive at
\begin{equation*}
	\p_x\widehat{U}(p,L) = \widehat{g}(p)\left[\widetilde{C}_0 + \sum_{i=0}^n \widetilde{C}_{i+1} p^{\alpha_i} + O\left(p^{2\min_i\{\alpha_i\}}\right)\right],
\end{equation*}
where the constants $\widetilde{C}_i$ are given by
\begin{equation}\label{eq:C0C1C2tild:def}
	\left\{\begin{aligned}
	\widetilde{C}_0&:={v_n'(L)}/{X_1^n},\\
	\widetilde{C}_{i+1}&:=(u,\overline{u})_{L^2_\rho(x_i,x_{i+1})}/\sigma(L)\quad\mbox{for }i=0,1,\dots,n.
\end{aligned}\right.
\end{equation}

The next result completes the proof of Theorem \ref{thm:asympt:1stord} on $\partial_x\widehat{U}(p,L)$.
\begin{lemma}\label{lem:X-repre-baru}
The function $\overline{u}$ defined in \eqnref{ubar:def:exp} belongs to $C^{1,1}[0,L]$ and satisfies
\beq\label{eq:ivp:ubar:C0dep}
    \left\{\begin{aligned}
	-(\sigma(x)\overline{u}'(x))' + q(x) \overline{u}(x)&=0,\quad\mbox{in }(0,L),\\
      \overline{u}(0)=0,\quad \overline{u}'(0)&=\frac{\sigma(L)}{\sigma(0)}\widetilde{C}_0.
		\end{aligned}\right.
\eeq   
\end{lemma}
\begin{proof}
Fix any $i=1,\dots,n$. Then we have
$$\lim_{s\to0^+}\overline{u}(x_i+s)=\frac{X_{1}^{i-1}}{X_1^n}\frac{v_n'(L)}{v_{i-1}'(x_i)}\frac{\sigma(L)}{\sigma(x_i)}=\lim_{s\to0^+}\overline{u}(x_i-s).$$ Thus, $\overline{u}$ is continuous at $x_i$. From the definition $X_1^i=c_i^*X_1^{i-1}+d_{i-1}^*X_1^{i-2}$, cf. \eqnref{eq:def:X}, we have
\begin{align*}
    \lim_{s\to0^+}\p_x \overline{u}(x_i + s)
    =&\frac{X_1^{i-1}}{X_1^n}\frac{v_n'(L)}{v_{i-1}'(x_{i})}\frac{\sigma(L)}{\sigma(x_{i})}v_i'(x_i)\\
     &+\frac{c_{i}^* X_1^{i-1}+d_{i-1}^* X_1^{i-2}}{X_1^n}\frac{v_n'(L)}{v_{i}'(x_{i+1})}\frac{\sigma(L)}{\sigma(x_{i+1})}w_i'(x_i)\\
    =&\frac{X_1^{i-1}}{X_1^n}\left(\frac{v_n'(L)}{v_{i-1}'(x_{i})}\frac{\sigma(L)}{\sigma(x_{i})}v_i'(x_i) + c_i^*\frac{v_n'(L)}{v_{i}'(x_{i+1})}\frac{\sigma(L)}{\sigma(x_{i+1})}w_i'(x_i)\right)\\
    &+\frac{X_1^{i-2}}{X_1^n}\left(d_{i-1}^*\frac{v_n'(L)}{v_{i}'(x_{i+1})}\frac{\sigma(L)}{\sigma(x_{i+1})}w_i'(x_i)\right).
\end{align*}
We analyze separately the two terms in brackets, denoted by ${\rm I}$ and ${\rm II}$. The definition of $c_i^*$ and Lemma \ref{lemma:Meq1} imply
\begin{align*}
  {\rm I}  =&\frac{v_n'(L)}{v_{i-1}'(x_{i})}\frac{\sigma(L)}{\sigma(x_{i})}\left(v_i'(x_i)+(v_i'(x_i)-w_{i-1}'(x_i))\frac{\sigma(x_i)w_i'(x_i)}{\sigma(x_{i+1})v_i'(x_{i+1})}\right)
    =\frac{v_n'(L)}{v_{i-1}'(x_i)}\frac{\sigma(L)}{\sigma(x_i)}w_{i-1}'(x_i).
\end{align*}
Likewise, using the definition of $d_i^*$ and Lemma \ref{lemma:Meq1}, we have
\begin{align*}
    {\rm II}=&\frac{w_{i-1}'(x_{i-1})}{v_{i-2}'(x_{i-1})}\frac{v_n'(L)}{v_{i}'(x_{i+1})}\frac{\sigma(L)}{\sigma(x_{i+1})}w_i'(x_i)\left(\frac{\sigma(x_i)v_{i-1}'(x_i)\sigma(x_{i+1})v_{i}'(x_{i+1})}{\sigma(x_{i-1})w_{i-1}'(x_{i-1})\sigma(x_{i})w_{i}'(x_{i})}\right)\\
    =&\frac{v_n'(L)}{v_{i-2}'(x_{i-1})}\frac{\sigma(L)}{\sigma(x_{i-1})}v_{i-1}'(x_i).
\end{align*}
Thus,  $\overline{u}'$ is continuous at $x_i$. 
Moreover,  $\overline{u}|_{[x_i,x_{i+1}]}$ belongs to $C^{1,1}[x_i,x_{i+1}]$, $i=0,1,\dots,n$. Then the gluing property for Lipschitz continuity (cf. the proof of Lemma \ref{lemma:u:ivp:C0dep}) implies $\overline{u}\in C^{1,1}[0,L]$. The identities in \eqnref{eq:ivp:ubar:C0dep} then follow from the definitions of $\overline{u}$ and $\widetilde{C}_0$.
\end{proof}

\begin{remark}
For the solution $U(t,x)$ of problem \eqnref{eq:withcoeff:ibvp}, the function $U(t,L-x)$ satisfies the governing equation of problem \eqnref{eq:withcoeff:ibvp} with $\alpha(L-x)$, $\rho(L-x)$, $\sigma(L-x)$ and $q(L-x)$ instead of $\alpha$, $\rho$, $\sigma$ and $q$, so the result below applies also to the case with the excitation $g$ at $x=L$ but the Neumann data at $x=0$.
\end{remark}
		
\section{Uniqueness of variable fractional order} \label{thmsec:Unique:invpb}
Now we prove Theorems \ref{thm:Unique:invpb} and \ref{thm:uniqueness:unknowncoeff} using Theorem \ref{thm:asympt:1stord}.

\subsection{Proof of Theorem \ref{thm:Unique:invpb}}

We prove the two cases of the theorem separately.

\medskip
\noindent\textbf{Proof of Case (i)}. By Theorem \ref{thm:weak}, $U^j$ is analytic in $t\in (0,\infty)$ as a map valued in $H^2(0,L^j)$. By the unique analytic continuation of $t\mapsto\p_x U^j(t,0)$, we have $\p_x U^1(t,0)=\p_x U^2(t,0)$ for all $t>0$, and thus,
$\p_x \widehat{U}^1(p,0)=\p_x \widehat{U}^2(p,0)$ for all $ p>0.$
To prove $\alpha^1=\alpha^2$,  we write, for $j=1,2$,
\begin{equation*} \alpha^j=\alpha_{0}^j+\sum_{i=1}^{n^j}(\alpha_{i}^j-\alpha_{0}^j)\chi_{(x_{i}^j,x_{i+1}^j)}\quad \mbox{and}\quad  L^j=x_{n^j+1}^j.
\end{equation*}
Without loss of generality, we may assume that $n^1\le n^2$ and for $j=1,2$,
$\alpha_{k}^j<\alpha_{l}^j$ for all $k,l\in\{0,1,\dots,n^j\}$ satisfying $k<l.$
From Theorem \ref{thm:asympt:1stord}, we have
\begin{equation}\label{eq:compare:cn}
	C_{0}^1 + \sum_{i=0}^{n^1} C_{i+1}^1 p^{\alpha_{i}^1} = C_0^2 + \sum_{i=0}^{n^2} C_{i+1}^2 p^{\alpha_{i}^2} + O\left(p^{\min\{2\alpha_{0}^1,2\alpha_{0}^2\}}\right),\quad \mbox{as } p\to0,
\end{equation}
where $C_{i+1}^j$ is defined in \eqnref{eq:C:intform} with $u^j$ satisfying	\begin{equation}\label{eq:initial:lappq}
	\left\{\begin{aligned}
			(\sigma(x)(u^j)'(x))' + q(x) u^j(x)&=0,\quad\mbox{in }(0,L^j),\\
			u^j(0)=1,\quad (u^j)'(0)&=C_{0}^j.
		\end{aligned}\right.
\end{equation}
It follows from \eqnref{eq:compare:cn} that $C_{0}^1=C_0^2=:C_0$, and thus $u^1=u^2=:u^0$. Meanwhile, the definition \eqnref{u:def:exp} and Lemma \ref{lem:X:positive} imply that $u^0$ is strictly positive in $[0,L^j)$ and $u^0(L^j)=0$ for both $j=1,2$, which gives $L^1=L^2$.  The positivity of $u^0$ and the expression \eqnref{eq:C:intform} give $C_{i+1}^j<0$ for $j=1,2$ and $i=0,1,\dots,n^j$.
From the identity \eqnref{eq:compare:cn}, we obtain $C_{1}^1=C_{1}^2$ and $\alpha_{0}^1=\alpha_{0}^2$.
Using condition \eqnref{eq:alpha:condition} and \eqnref{eq:compare:cn}, we also get $C_{i+1}^1=C_{i+1}^2$, $\alpha_{i}^1=\alpha_{i}^2$ for $i=1,\dots,n^1$ and $n^1=n^2=:n$.
Next, we rewrite the expression of $C_{i+1}^j$ in \eqnref{eq:C:intform} for $i=0,1,\dots,n$:
$$C_{i+1}^j	= - \frac{1}{\sigma(0)}\int_{x_{i}^j}^{x_{i+1}^j} |u^0(x)|^2\rho(x)\,{\rm d}x.$$
Now for $i=0,1,\dots,n-1$, $C_{i+1}^j$ is a function of $x_{i}^j$ and $x_{i+1}^j$ satisfying
\begin{equation} \label{eq:C:monotone}
\frac{\p C_{i+1}^j}{\p x_{i+1}^j}= -\frac{1}{\sigma(0)}|u_0(x_{i+1}^j)|^2\rho(x_{i+1}^j)<0, \quad\mbox{for all }x_{i+1}^j>x_{i}^j.
\end{equation}
Thus, for each $i=0,1,\dots,n-1$, the identity $C_{i+1}^1=C_{i+1}^2$ implies
$x_{i+1}^1=x_{i+1}^2$, i.e., $\alpha^1=\alpha^2$.

\medskip
\noindent\textbf{Proof of Case (ii).}
Let $V^j$ with $j=1,2$ be defined by
$V^j(t,x):=U^j(t,L^j-x)$ for all $(t,x)\in(0,\infty)\times(0,L^j).$
Then $V^j$ is the solution to problem \eqnref{eq:withcoeff:ibvp} with $\alpha^j(L^j-x)$, $\rho(L^j-x)$, $\sigma(L^j-x)$ and $q(L^j-x)$ in place of $\alpha$, $\rho$, $\sigma$ and $q$.
By applying Theorem \ref{thm:asympt:1stord} to $V^j$ for $j=1,2$, we have
$$\p_x\widehat{U}^j(p,0)=-\p_x\widehat{V}^j(p,L^j)=-\widehat{g}(p)\left( \widetilde{C}_{0}^j + \sum_{i=0}^{n^j} \widetilde{C}_{i+1}^j p^{\alpha_{i}^j} + O\left(p^{\min\{2\alpha_{i}^j\,:\,0\le i\le n^j\}}\right)\right),$$
where $\widetilde{C}_{i+1}^j$ has an integral expression, cf. \eqnref{eq:C:intform}, in terms of the unique $C^{1,1}$ solutions $u^j$ and $\overline{u}^j$ of the following problems
		\begin{equation*}
			\left\{\begin{aligned}
				-(\sigma (L^j-x)(u^j)'(x))' + q(L^j-x) u^j(x)&=0,\quad\mbox{in }(0,L^j),\\
				u^j(0)=1,\quad u^j(L^j)&=0,
		\end{aligned}\right.
		\end{equation*}
		\begin{equation*}
			\left\{\begin{aligned}
				-(\sigma (L^j-x)(\overline{u}^j)'(x))' + q(L^j-x) \overline{u}^j(x)&=0,\quad\mbox{in }(0,L^j),\\
				\overline{u}^j(0)=0,\quad \overline{u}^j(L^j)&=1.
		\end{aligned}\right.
		\end{equation*}
    The function $u^j$ is strictly decreasing in $(0,L^j)$ (since $u^j$ is a fundamental solution in $(0,L^j)$, cf. the proof of Lemma \ref{lem:v:decreasing}) and satisfies $(u^j)'(L^j)=\widetilde{C}_{0}^j$, so the identity $\widetilde{C}_{0}^1=\widetilde{C}_{0}^2$ gives $u^1(L^1-y)=u^2(L^2-y)$ for all $y>0$ so that
    $$L^1=\min\{y>0\,:\,u^1(L^1-y)=1\}=\min\{y>0\,:\,u^2(L^2-y)=1\}=L^2.$$
Thus, we deduce $\overline{u}^1=\overline{u}^2$. Now we proceed as in the proof of case (i), using the monotonicity of $\alpha^j$ and condition \eqnref{eq:alpha:condition}, and identify $n^1=n^2=:n$ and
$\widetilde{C}_{i+1}^1=\widetilde{C}_{i+1}^2$ and $\alpha_{i}^1=\alpha_{i}^2$ for $ i=1,\dots,n.$ Since we have $u^1=u^2>0$ and $\overline{u}^1=\overline{u}^2>0$, the identity  $\widetilde{C}_{i+1}^1=\widetilde{C}_{i+1}^2$ implies
$x_{i+1}^1=x_{i+1}^2$ for each $i=0,1,\dots,n-1$.
Thus we conclude $\alpha^1=\alpha^2$.
 
\subsection{Proof of Theorem  \ref{thm:uniqueness:unknowncoeff}}

The proof of Theorem \ref{thm:uniqueness:unknowncoeff} mainly relies on Theorem \ref{thm:asympt:1stord}. 

\begin{proof}[Proof of Theorem \ref{thm:uniqueness:unknowncoeff}]
The proof proceeds as that of Theorem \ref{thm:Unique:invpb}. By the unique analytic continuation of $t\mapsto\p_x U^j(t,0)$,
\begin{equation}\label{eq:neumann:unknowncoeff}
\p_x \widehat{U}^1(p,0)=\p_x \widehat{U}^2(p,0)\quad\mbox{for all }p>0.
\end{equation}
We will prove the assertion \eqnref{eq:alpha:identicalset} for the cases $(a^j,b^j)=(0,L^j)$ and $(a^j,b^j)=(L^j,0)$ separately.
            
\medskip
\noindent{\bf Case 1. ($(a^j,b^j)=(0,L^j)$)}
Suppose on the contrary that the assertion \eqnref{eq:alpha:identicalset} does not hold. We may assume that
\beq\label{eq:contrad:beta:hyp}
\beta\in\{\alpha^1(x)\,:\,x\in(0,L^1)\}\quad\mbox{but}\quad \beta\not\in\{\alpha^2(x)\,:\,x\in(0,L^2)\}.
\eeq
It follows from condition \eqnref{eq:g:poly} that
\begin{equation}\label{eqn:hat-g} \widehat{g^j}(p)=\sum_{k=2}^{N^j} g_{k}^j\frac{k!}{p^{k+1}}= \frac{N^j! g_{N_j}^j}{p^{N^j+1}}\left(1+O(p)\right)\quad\mbox{as }p\to0.
\end{equation}
From \eqnref{eq:neumann:unknowncoeff}, as $p\to0$, there holds
\begin{align*}
    &\frac{N^1! g_{N^1}^1}{p^{N^1+1}}\bigg(C_{0}^1 + \sum_{i=0}^{n^1} C_{i+1}^1 p^{\alpha_{i}^1}\bigg)\left(1+O\left(p+p^{\min\{2\alpha_{i}^1\,:\,0\le i\le n^1\}}\right)\right) \\
    =& \frac{N^2! g_{N^2}^2}{p^{N^2+1}}\bigg(C_{0}^2 + \sum_{i=0}^{n^2} C_{i+1}^2 p^{\alpha_{i}^2}\bigg)\left(1+O\left(p+p^{\min\{2\alpha_{i}^2\,:\,0\le i\le n^2\}}\right)\right).
            \end{align*}
Since the function $\alpha^j$ is valued in $(0,1)$ for both $j=1,2$, we have $N^1=N^2=N$.
Since $\alpha=\alpha^j$ satisfies condition \eqnref{eq:alpha:condition}, there holds
\beq\label{eq:same:coeff:thm1p4}
\sum_{i:\,\alpha_{i}^1=\beta} g_{N}^1C_{i+1}^1=\sum_{i:\,\alpha_{i}^2=\beta}  g_{N}^2C_{i+1}^2.
\eeq
By Theorem \ref{thm:asympt:1stord}, we have $C_{i+1}^j<0$ for $j=1,2$ and $i=0,1,\dots,n^j$.
By the hypothesis \eqnref{eq:contrad:beta:hyp}, the left-hand side of \eqnref{eq:same:coeff:thm1p4} is nonzero, but the right-hand side is zero, which leads to a contradiction.

\smallskip\smallskip
\noindent{\bf Case 2. ($(a^j,b^j)=(L^j,0)$)}
The proof in this case is identical with that for the case with $(a^j,b^j)=(0,L^j)$ and $\p_x \widehat{U}^1(p,L^1)=\p_x \widehat{U}^2(p,L^2)$ for all $p>0$, but using $\widetilde{C}_{i+1}^j>0$ instead of $C_{i+1}^j<0$.
\end{proof}


The last result follows from Theorem \ref{thm:uniqueness:unknowncoeff}, dealing with a constant $\rho$. It gives the unique determination of the constant $\rho$ and the length $L$ of the interval $(0,L)$ without the monotonicity of the variable order $\alpha$.
\begin{cor}
Let $L^j$, $g$, $\sigma$ and $q$ satisfy the assumptions of Theorem \ref{thm:Unique:invpb}. Let $\alpha^j:(0,L^j)\to(0,1)$ be piecewise constant and satisfy \eqref{eq:alpha:condition} and let $\rho^j>0$ be constant for $j=1,2$. Let $U^j$ be defined as in Theorem \ref{thm:Unique:invpb} with $\rho=\rho^j$. If $\p_x U^1(t_k,0)=\p_x U^2(t_k,0)$ for some sequence $\{t_k\}_{k=1}^\infty$ of distinct numbers converging to a positive number, then we have $L^1=L^2$ and $\rho^1=\rho^2$.
\end{cor}

\begin{proof}
We first prove the assertion for the case $U^j(t,0)=g(t)$ and $U^j(t,L^j)=0$ for all $t>0$. 
By Theorem \ref{thm:asympt:1stord}, we have $C_{0}^1=C_{0}^2$.
Then we have $u^1=u^2$ so that $L^1=\min\{x>0\,:\,u^1(x)=0\} =\min\{x>0\,:\,u^2(x)=0\} =L^2$.
Also, from Theorem \ref{thm:asympt:1stord}, \eqnref{eq:alpha:identicalset} and condition \eqnref{eq:alpha:condition}, we have
$$\frac{1}{\sigma(0)}\int_0^{L^1} |u^1(x)|^2\rho^1\,{\rm d}x =\sum_{i=1}^{n^1} C_{i+1}^1=\sum_{i=1}^{n^2} C_{i+1}^2=\frac{1}{\sigma(0)}\int_0^{L^2} |u^2(x)|^2\rho^2\,{\rm d}x.$$
Since $L^1=L^2>0$ and $u^1=u^2>0$, there holds $\rho^1=\rho^2$.
Next, we prove the assertion for the case when $U^j(t,0)=0$ and $U^j(t,L^j)=g(t)$ for all $t>0$.
Similar to the proof of case (ii) of Theorem \ref{thm:Unique:invpb}, we define $V^j(t,x):=U^j(t,L^j-x)$ and $\widetilde{C}_{i}^j$ for $j=1,2$.
Similarly, from $\widetilde{C}_{0}^1=\widetilde{C}_{0}^2$, we have $u^1(L^1-y)=u^2(L^2-y)$ for all $y>0$, which gives $L^1=\min\{y>0\,:\,u^1(L^1-y)=0\} =\min\{y>0\,:\,u^2(L^1-y)=0\} =L^2$.
We therefore have $\overline{u}^1=\overline{u}^2$.
From Theorem \ref{thm:asympt:1stord}, \eqnref{eq:alpha:identicalset} and condition \eqnref{eq:alpha:condition}, we have
$$\frac{1}{\sigma(L^1)}\int_0^{L^1} u^1(x)\overline{u}^1(x)\rho^1\,{\rm d}x =\sum_{i=1}^{n^1} \widetilde{C}_{i+1}^1=\sum_{i=1}^{n^2} \widetilde{C}_{i+1}^2=\frac{1}{\sigma(L^2)}\int_0^{L^2} u^2(x)\overline{u}^2(x)\rho^2\,{\rm d}x.$$
Since $L^1=L^2>0$, $u^1=u^2>0$ and $\overline{u}^1=\overline{u}^2>0$, there holds $\rho^1=\rho^2$.
\end{proof}

\section{Conclusion}
In this work we have investigated of the identifiability issue of the spatially variable-order in one-dimensional subdiffusion, and proved the uniqueness of a monotone piecewise constant variable order from one boundary flux observation. The main tool in the analysis is a new first-order asymptotic expansion of the boundary data in the Laplace domain. The work contributes to the study on inverse problems for variable order models, which remains fairly scarce so far. Even for the setting of a piecewise constant variable order in the one-dimensional case, there are still several outstanding questions to be further studied, including stability of the recovery, simultaneous recovery of the order and other parameters, and efficient algorithms for recovery. The analysis in this work relies crucially on a gluing strategy (i.e., matching the continuity conditions at intersection points of the subintervals) and local Sturm-Liouville problems on the subintervals, and thus it does not extend directly to general multi-dimensional cases, for which we expect that a completely new strategy is needed. It is of great interest to investigate recovering the support of the variable-order $\alpha(x)$ from one lateral flux measurement in the multi-dimensional case either mathematically or numerically.

\bibliographystyle{siam}
\bibliography{frac}
\end{document}